\documentclass[11pt]{article}
\usepackage{amssymb}
\usepackage{latexsym}
\usepackage{amsfonts}
\usepackage{amsmath}
\newtheorem{thrm}{Theorem}[section]
\newtheorem{lemma}[thrm]{Lemma}
\newtheorem{prop}[thrm]{Proposition}

\newtheorem{cor}[thrm]{Corollary}
\newtheorem{remark}[thrm]{Remark}
\numberwithin{equation}{section}

\usepackage[dvips]{color}

\def\P{\mathbb{P} }
\def\R{\mathbb{R} }
\def\N{\mathbb{N} }

\def\bP{{\bf P} }

\def\cL{{\cal L} }

\begin{document}
\allowdisplaybreaks
\begin{doublespace}

\title{\Large\bf Central Limit Theorems for Super-OU Processes}
\author{ \bf  Yan-Xia Ren\footnote{The research of this author is supported by NSFC (Grant No.  10871103 and 10971003) and Specialized Research Fund for the
Doctoral Program of Higher Education.\hspace{1mm} } \hspace{1mm}\hspace{1mm}
Renming Song\thanks{Research supported in part by a grant from the Simons
Foundation (208236).} \hspace{1mm}\hspace{1mm} and \hspace{1mm}\hspace{1mm}
Rui Zhang\footnote{Supported by the China Scholarship Council}
\hspace{1mm} }
\date{}
\maketitle
\begin{abstract}

In this paper we study supercritical super-OU processes with general branching mechanisms satisfying a
second moment condition. We establish central limit theorems for the super-OU processes. In the small and crtical
branching rate cases, our central limit theorems sharpen the corresponding results in the recent preprint of
Milos in that the limit normal random variables in our central limit theorems are non-degenerate. Our
central limit theorems in the large branching rate case are completely new.
The main tool of the paper is the so called ``backbone decomposition" of superprocesses.
\end{abstract}


\medskip


\bigskip

\baselineskip=6.0mm

\section{Introduction}
\subsection{Model}

Throughout this paper, $d\ge1$ is a positive integer and $b$ is a positive number.
We use $\xi=\{\xi_t: t\ge0\}$
to denote an Ornstein-Uhlenbeck process (OU process, for short)
on $\mathbb{R}^d$, that is, a diffusion process with infinitesimal generator
\begin{equation}\label{2}
  L:=\frac{1}{2}\sigma^2\bigtriangleup-b x\cdot\bigtriangledown.
\end{equation}
For any $x\in \mathbb{R}^d$, we use $\Pi_x$ to denote the law of $\xi$ starting from $x$.
The semigroup of $\xi$ will be denoted by $\{T_t:t\geq0\}$.

Consider a branching mechanism of the form
\begin{equation}
\psi(\lambda)=-\alpha\lambda+\beta\lambda^2+\int_{(0,+\infty)}(e^{-\lambda x}-1+\lambda x)n(dx),
\quad \lambda> 0,
\end{equation}
where $\alpha=-\psi'(0+) >0$, $\beta\geq 0$, and $n$ is a measure on $(0,\infty)$ such that
\begin{equation}\label{cond-n}
\int_{(0,+\infty)} x^2n(dx)<+\infty.
\end{equation}

Let ${\cal M}_F(\mathbb{R}^d)$ be the space of finite  measures on $\mathbb{R}^d$.
In this paper we will always assume that $X=\{X_t: t\ge 0\}$ is a
super-Ornstein-Uhlenbeck process (super-OU process, for short) with underlying spatial
motion $\xi$ and branching mechanism $\psi$. We will sometimes call $X$ a
$(\xi, \psi)$-superprocess.
The existence of such superprocesses is well-known, see, for instance, \cite{E.B.}.
$X$ is a Markov branching branching process taking values in ${\cal M}_F(\mathbb{R}^d)$.
For any $\mu \in \mathcal{M}_F(\mathbb{R}^d)$, we denote the
law of $X$ with initial configuration $\mu$ by $\P_\mu$.
The total mass of the process $X$ is a continuous-state branching process
with branching mechanism $\psi$.
The assumption \eqref{cond-n} implies that the total mass process of $X$
does not explode. Since we always assume that $\alpha>0$, $X$ is a supercritical
superprocess.

Let $\mathcal{B}$$^+_b(\mathbb{R}^d)$ be the space of positive, bounded
measurable functions on $\mathbb{R}^d$. As usual, $\langle f,\mu\rangle
:=\int f(x)\mu(dx)$ and $\|\mu\|:=\langle 1,\mu\rangle$. Then for every
$f\in \mathcal{B}^+_b(\mathbb{R}^d)$ and $\mu \in \mathcal{M}_F(\mathbb{R}^d)$,
\begin{equation}
  -\log \P_\mu\left(e^{-\langle f,X_t\rangle}\right)=\langle u_f(\cdot,t),\mu\rangle,
\end{equation}
where $u_f(x,t)$ is the unique positive solution to the equation
\begin{equation}
  u_f(x,t)+\Pi_x\int_0^t\psi(u_f(\xi_s,t-s))ds=\Pi_x f(\xi_t).
\end{equation}

In addition, we assume that $\psi(\infty)=\infty$
which implies that the probability of the extinction event
$\mathcal{E}:=\{\lim_{t\to \infty} \|X_t\|=0\}$ in strictly in $(0, 1)$,
see for example the summary at then end of \cite[Section 10.2.2]{A.E.}.
Since $\psi$ is convex with $\psi(0)=0,\psi(\infty)=\infty$
and $\psi'(0+)<0$, $\psi$ has exactly two roots in $[0,\infty)$;
let $\lambda^*$ be the larger one. We have
\begin{equation*}
  \P_\mu(\lim_{t\to\infty}\|X_t\|=0)=e^{-\lambda^*\|\mu\|}.
\end{equation*}
Using the expectation formula of $\|X_t\|$ and the Markov property of $X$, it is not hard to prove that (see Lemma \ref{lem:1.2} for a proof), under $\P_{\mu}$,
the process $W_t=e^{-\alpha t}\|X_t\|$ is a positive martingale. Therefore it converges:
\begin{equation}\label{3.1}
  W_t \to W_\infty,\quad\P_{\mu}\mbox{-a.s.} \quad \mbox{ as }t\to \infty.
\end{equation}
Using the assumption \eqref{cond-n}
we can show that, as $t\to \infty$, $W_t$ also converges in $L^2(\P_\mu)$, so $W_\infty$ is non-degenerate and
 the second moment is finite. Moreover, we have $\P_{\mu}(W_\infty)=\|\mu\|$ and
 $\{W_\infty=0\}=\mathcal{E}$.

The purpose of this paper is to establish central limit theorems for the super-OU process.
More precisely, we want to find $A_t$ and $C_t$, for suitable test functions $f$, such that
$C_t(\langle f,X_t\rangle -A_t)$ converges to  some non-degenerate random variable as $t\to \infty$.
It turns out that $C_t$ is determined by the second moment of $\langle f,X_t\rangle$
which depends on the sign of $\alpha-2\gamma(f)b$, where $\gamma(f)$ is a quantity to
be defined later.

There are many papers studying laws of large numbers for branching processes and superprocesses.
For example, see \cite{AH76, AH77, EHK} for branching processes,
and ~\cite{EW, E09, LRS2} for superprocesses.
For super-OU processes with binary branching mechanism, the following weak law of large numbers
was proved in \cite{EW}:
\begin{equation}\label{1.1}
  e^{-\alpha t}\langle f,X_t\rangle\to \langle f,\varphi\rangle W_\infty, \quad \mbox{in probability}
\end{equation}
where $f\in C^+_c(\mathbb{R}^d)$.
When $\langle f,\varphi\rangle=0$, it is natural to consider central limit theorems for $\langle f,X_t\rangle$,
that is,
to find a normalization $C_t$ so that $C_t\langle f,X_t\rangle$ converges to a non-degenerate
Gaussian random variable as $t\to\infty$.
For branching processes, there are already papers dealing with central limit theorems.
In 1966, Kesten and Stigum \cite{KS66} gave a central limit theorem for multidimensional
Galton-Watson processes by using the Jordan canonical form of the expectation matrix $M$.
Then in \cite{Ath69a, Ath69, Ath71},
Athreya proved central limit theorems
for multi-type continuous time Markov branching processes; the main tools used
in \cite{Ath69a, Ath69, Ath71}
are also the Jordan canonical form and the eigenvectors of the matrix $M_t$, the mean matrix at time $t$.
Recently, central limit theorems for branching OU particle systems
and for super-OU processes were established in \cite{RP} and \cite{Mi} respectively.
However, the limiting normal random variables in the central limit theorems in \cite{RP, Mi}
(see \cite[Theorems 3.2 and 3.6]{RP} and \cite[Theorems 3 and 4]{Mi}) may be degenerate (i.\/e., equal to
zero), so the central limit theorems in \cite{RP, Mi} are not completely satisfactory.

In this paper, we sharpen the results of \cite{Mi} and establish central limit theorems for super-OU processes
which are more satisfactory in the sense that
the limiting normal random variables in our results are non-degenerate.
The setup of this paper is more general than that of \cite{Mi} in the sense that we allow a general branching mechanism
as opposed to the binary branching mechanism in \cite{Mi}.
The only assumption on the branching mechanism is the second momemnt condition
\eqref{cond-n}, which is necessary for central limit theorems.

We mention that we are following Athreya's argument for multitype (finite type) branching processes,
also called  multidimensional Galton-Watson processes,
and show that Athreya's ideas for finite dimensional branching processes also work for super-OU processes,
which can be regarded as an infinite dimensional branching process.
The main tool of this paper is, similar to that of \cite{Mi},
also the backbone decomposition of supercritical superprocesses, see~\cite{BKS}.
The main idea of the backbone decomposition is that a supercritical
super-OU process can be constructed from a branching OU process (known as the backbone),
in which particles live forever (known as immortal particles).
After dressing the backbone with subcritical super-OU processes,
we get a measure-valued Markov
process,
which gives a version of the super-OU process. Since subcritical super-OU process will die
out in finite time, we can imagine that the
limit behavior of super-OU process is determined by the backbone branching OU process.
In this paper we prove that these intuitive ideas work well.
For the precise backbone decomposition, see Section 2.1.

Using a similar argument, we can also sharpen results of \cite{RP} and establish central limit theorems for
branching OU particle systems which are more satisfactory in the sense that limiting normal random variables
are non-degenerate.

\subsection{Eigenfunctions of OU processes}

Recall that $\{T_t,t\geq0\}$ is the semigroup of the OU process $\xi$.
It is well known that $\xi$ has an invariant density
\begin{equation}\label{2.1}
\varphi(x)=\left(\frac{b}{\pi\sigma^2}\right)^{d/2}\exp \left(-\frac{b}{\sigma^2}\|x\|^2\right).
\end{equation}
Let $L^2(\varphi):=\{h: \int_{\mathbb{R}^d}|h(x)|^2\varphi(x)dx<\infty\}$. For $h_1,h_2\in L^2(\varphi)$, we define
\begin{equation*}
\langle h_1, h_2\rangle_\varphi:=\int_{\mathbb{R}^d} h_1(x)h_2(x)\varphi(x)dx.
\end{equation*}
In this subsection, we recall some results on the spectrum in $L^2(\varphi)$ of the operator $L$ defined in \eqref{2},
more details can be found in~\cite{MP}.
For $p=(p_1, p_2,\cdots, p_d)\in \mathbb{Z}_+^d$, let $|p|=\sum_{j=1}^{d}p_j$ and $p!=\prod_{j=1}^{d}p_j!$. Recall the Hermite polynomials $\{H_p(x), ~p\in \mathbb{Z}_+^d\}$:
\begin{equation}\label{2.4}
  H_p(x)=(-1)^{\mid p\mid}e^{\| x\|^2}\frac{\partial}{\partial x_1^{p_1}\cdots\partial x_d^{p_d} }(e^{-\| x\|^2}).
\end{equation}
The eigenvalues of $L$
are $\{-mb,m=0,1,2...\}$ and the corresponding eigenspaces $A_m$ are given by
\begin{equation}\label{2.3}
  A_m:= {\rm Span }\{\phi_p,~|p|=m\},
\end{equation}
where
\begin{equation*}
\phi_p(x)=\frac{1}{\sqrt{p!2^{|p|}}}H_p\left(\frac{\sqrt{b}}\sigma x\right).
\end{equation*}
In particular, $\phi_{0,0,\cdots,0}(x)=1$, $\phi_{e_i}(x)=\frac{\sqrt{2b}}{\sigma}x_i$,
where $e_i$ stands for the unit vector in the $x_i$ direction.
The function $\phi_p$ is an eigenfunction of $L$ corresponding to the eigenvalue $-|p|b$ and therefore
\begin{equation}\label{2.18}
  T_t\phi_p(x)=e^{-|p|b t}\phi_p(x).
\end{equation}
Moreover, the eigenfunctions $\{\phi_p(x), ~p\in \mathbb{Z}_+^d\}$
form a complete orthonormal basis for $L^2(\varphi)$. Thus every $f\in L^2(\varphi)$ admits the
following $L^2(\varphi)$ expansion:
\begin{equation}\label{2.5}
  f(x)=\sum_{m=0}^{\infty}\sum_{|p|=m}a_p\phi_p(x),
\end{equation}
where $a_p=\langle f,\phi_p\rangle_\varphi$. Define
\begin{eqnarray}\label{2.19}
  \gamma(f):=\inf\{n\geq 0: \mbox{ there exists } p\in \mathbb{Z}_+^d \mbox{ with }
  |p|=n \mbox{ such that }
  a_p\neq 0\},
\end{eqnarray}
where we use the usual convention $\inf\varnothing=\infty$. Define
$$
  f_{(s)}(x) = \sum_{\gamma(f)\le m<\alpha/(2b)}\sum_{|p|=m}a_p\phi_p(x),
  \quad f_{(c)}(x)=\sum_{m=\alpha/(2b)}\sum_{|p|=m}a_p\phi_p(x),
 $$ and $$
  f_{(l)}(x)=f(x)-  f_{(s)}(x)-f_{(c)}(x)=\sum_{m>\alpha/(2b)}^\infty\sum_{|p|=m}a_p\phi_p(x).
$$
In this paper we will use $\mathcal{P}$ to denote the function class
\begin{equation*}
  \mathcal{P}:=\{f\in C(\mathbb{R}^d): \mbox{ there exists } ~k\in\mathbb{N} \mbox{  such that }|f(x)|/\|x\|^k\to 0
   \mbox{ as }\|x\|\to\infty \}.
\end{equation*}
We easily see that $\mathcal{P}\subset L^2(\varphi)$ and for $f\in \mathcal{P}$, there exists $k\in \mathbb{N}$,
\begin{equation*}
  |f(x)|\lesssim 1+\|x\|^k,
\end{equation*}
where we used the following notation: for two positive functions $f$ and $g$,
$f(x)\lesssim g(x)$ means that there exists a constant $c>0$ such that $f(x)\leq c g(x)$.

\subsection{Main results for super-OU processses}

In this subsection we give the main results of this paper.
The proofs will be given in the later sections.
In the remainder of this paper,
whenever we deal with
an initial configuration $\mu\in {\cal M}_F(\mathbb{R}^d)$, we
are implicitly assuming that it has compact support.

\subsubsection{Large branching rate: $\alpha > 2b\gamma(f)$}

For each $p\in \mathbb{Z}^d_+$, we define
\begin{equation*}
  H_t^{p}:=e^{-(\alpha-|p|b)t}\langle\phi_p, X_t\rangle,\quad t\geq 0.
\end{equation*}
Then one can show (see Lemma \ref{lem:1.2} below) that, if $\alpha>2|p|b$,
then under $\P_{\mu}$,
$H_t^p$ is a martingale bounded in $L^2(\P_{\mu})$, and thus the limit
$ H_\infty^p:=\lim_{t\to \infty}H_t^p$ exists $\P_{\mu}$-a.s. and in $L^2(\P_{\mu})$.

\begin{thrm}\label{The:1.3}
If  $f\in \mathcal{P}$ satisfies $\alpha > 2\gamma(f)b$, then as $t\to\infty$,
\begin{equation*}
e^{-(\alpha-\gamma(f)b)t)}\langle f,X_t\rangle\to \sum_{|p|=\gamma(f)}a_p H_\infty^p,\quad \mbox{in }L^2(\P_{\mu}).
\end{equation*}
\end{thrm}

\begin{remark}\label{rem:large} When $\gamma(f)=0$, $H_t^0$ reduces to $W_t$, and thus $H_\infty^0=W_\infty$.
Therefore by Theorem~\ref{The:1.3} and the fact that $a_0=\langle f,\varphi\rangle$,
we get that, as $t\to\infty$,
\begin{equation*}
  e^{-\alpha t}\langle f,X_t\rangle\to \langle f,\varphi\rangle W_\infty, \quad \mbox{in } L^2(\P_{\mu}).
\end{equation*}
In particular, the convergence also holds in $\P_{\mu}$-probability,
so it implies the results in~\cite{EW} in the case of super-OU processes.
Moreover, by \eqref{3.1}, on $\mathcal{E}^c$, we have
\begin{equation*}
  \|X_t\|^{-1}\langle f,X_t\rangle\to \langle f,\varphi\rangle, \quad \mbox{ in } \P_{\mu}\mbox{-probability}.
\end{equation*}
\end{remark}

\subsubsection{Small branching rate: $\alpha<2\gamma(f)b$}

Let
\begin{equation}\label{e:sigmaf}
  \sigma_f^2:=A\int_0^{\infty} e^{\alpha s}\langle(T_sf)^2,\varphi\rangle\, ds,
\end{equation}
where
\begin{equation}\label{e:A}
A:=\psi^{(2)}(0+)=2\beta+\int_{(0,\infty)}x^2n(dx)<\infty.
\end{equation}
In the rest of this paper, $A$ will always stand for this constant.

\begin{thrm}\label{The:1.4}
If $f\in \mathcal{P}$ satisfies $\alpha < 2\gamma(f)b$, then $ \sigma_f^2<\infty$ and,
under $\P_{\mu}(\cdot\mid \mathcal{E}^c)$, it holds that
\begin{equation}\label{the:4.2}
  \left(e^{-\alpha t}\|X_t\|, ~\frac{\langle f , X_t\rangle}{\sqrt{\|X_t\|}} \right)\stackrel{d}{\rightarrow}(W^*,~G_1(f)), \quad t\to\infty,
\end{equation}
where $W^*$ has the same distribution as $W_\infty$ conditioned on $\mathcal{E}^c$
and $G_1(f)\sim \mathcal{N}(0,\sigma_f^2)$. Moreover, $W^*$ and $G_1(f)$ are independent.
\end{thrm}
\begin{remark} Using the theorem above, we get that if $\alpha<2\gamma(f)b$, then,
under $\mathbb{P}_\mu$, we have
\begin{equation*}
   e^{-\alpha t/2}\langle f,X_t\rangle\stackrel{d}{\rightarrow} G_1(f)\sqrt{W_\infty},
\end{equation*}
where $W_\infty$ and $G_1(f)$ are the same as in the theorem above.
\end{remark}

\subsubsection{The critical case: $\alpha = 2\gamma(f)b$}
Define
\begin{equation}\label{1.3}
  \rho_f^2:=A \sum_{|p|=\gamma(f)}(a_p)^2.
\end{equation}
\begin{thrm}\label{The:1.5}
If $f\in \mathcal{P}$ satisfies $\alpha = 2\gamma(f)b$, then, under $\P_{\mu}(\cdot\mid \mathcal{E}^c)$,
it holds that
\begin{equation*}
  \left(e^{-\alpha t}\|X_t\|, ~\frac{\langle f , X_t\rangle}{t^{1/2}\sqrt{\|X_t\|}} \right)\stackrel{d}{\rightarrow}(W^*,~G_2(f)), \quad t\to\infty,
\end{equation*}
 where $W^*$ has the same distribution as $W_\infty$ conditioned on $\mathcal{E}^c$,
$G_2(f)\sim \mathcal{N}(0, \rho_f^2)$. Moreover $W^*$ and $G_2(f)$ are independent.
\end{thrm}

\begin{remark} Using the theorem above, we get that if $\alpha=2\gamma(f)b$, then,
under $\mathbb{P}_\mu$, we have
\begin{equation*}
   t^{-1/2}e^{-\alpha t/2}\langle f,X_t\rangle\stackrel{d}{\rightarrow} G_2(f)\sqrt{W_\infty}, \quad t\to\infty,
\end{equation*}
where $W_\infty$ and $G_2(f)$ are the same as in the theorem above.
\end{remark}
\begin{remark}
Note that the limiting normal random variables in our Theorems \ref{The:1.4} and \ref{The:1.5}
are non-degenerate.
\end{remark}

\begin{remark}
The results of \cite{Mi} correspond to the case $\gamma(f)=1$ in the present paper.
For the small branching rate
case of \cite{Mi}, the $\sigma_f^2$ in (3.1) there should be (in the notation there)
\begin{equation*}
  \sigma_f^2=2\beta\int_0^\infty e^{-\alpha s}\langle \varphi,( \mathcal{P}^\alpha_s\widetilde{f}(\cdot))^2\rangle\,ds,
\end{equation*}
$\widetilde f(x)=f(x)-\langle f,\phi\rangle.$
It is easy to check that the sum of the last two parts of \cite[(3.1)]{Mi} is 0, that is
 \begin{equation*}
   \int_0^\infty\langle \varphi,(-2\beta( \mathcal{P}^{-\alpha}_s\widetilde{f}(\cdot))^2+4\alpha\beta u(\cdot,s))\rangle\,ds=0,
 \end{equation*}
 where $u(x,s)=\int_0^s( \mathcal{P}^{-\alpha}_{s-u}[( \mathcal{P}^{-\alpha}_u\widetilde{f}(\cdot))^2](x)\,du$.
 Furthermore, there is an extra factor $\beta/\alpha$ on the right side of \cite[(3.1)]{Mi} which should not be there.
  In the critical branching case of \cite{Mi}, there is also an extra factor $\beta/\alpha$ on the right side of \cite[(3.2)]{Mi} which should not be there. The correct form of (3.2) there should be (in the notation of \cite{Mi})
 \begin{equation*}
   \sigma_f^2=2\beta\int_{\mathcal{R}^d}(x\circ\langle grad( f),\varphi\rangle)^2\varphi(x)dx.
 \end{equation*}
 With these minor corrections, the results of \cite{Mi} coincide with our Theorems \ref{The:1.3}, \ref{The:1.4} and \ref{The:1.5} when $\gamma(f)=1$.
 \end{remark}

Combining Theorems \ref{The:1.3}, \ref{The:1.4} and \ref{The:1.5},
we have the following expansion  of $\langle f,X_t\rangle$: for any $f\in \mathcal{P}$,
 \begin{eqnarray}\label{expansion}
 \langle f, X_t\rangle&=&\sum_{\gamma(f)\le m<\frac{\alpha}{2b}}\sum_{|p|=m}a_pe^{-(\alpha-mb)t}\langle\phi_p, X_t\rangle\cdot e^{(\alpha-m)bt}\nonumber\\
&&+\sum_{|p|=\frac{\alpha}{2b}}a_pt^{-1/2}e^{-(\alpha/2) t}\langle \phi_p, X_t\rangle\cdot \sqrt{t}e^{(\alpha/2) t}
+\langle f_{(l)}, X_t\rangle\nonumber\\
&=&\sum_{\gamma(f)\le m<\frac{\alpha}{2b}}\sum_{|p|=m}a_pU_{p}(t)\cdot e^{(\alpha-m)bt}
+\sum_{|p|=\frac{\alpha}{2b}}a_pU_p(t)\cdot \sqrt{t}e^{\alpha t/2}
+\langle f_{(l)}, X_t\rangle,
\end{eqnarray}
where
$$U_p(t)=\left\{\begin{array}{ll}\displaystyle e^{-(\alpha-|p|b)t}\langle\phi_p, X_t\rangle,\quad &|p|<\frac{\alpha}{2b},\\
\displaystyle t^{-1/2}e^{-\alpha t/2}\langle \phi_p, X_t\rangle,\quad& |p|=\frac{\alpha}{2b}.
\end{array}\right.$$
Further, if $|p|<\frac{\alpha}{2b}$,
then $U_p(t)=H^{p}_t$ converges to $H^p_\infty$,
$\P_\mu$-a.s. and in $L^2(\P_\mu)$; if $|p|=\frac{\alpha}{2b}$, $U_p(t)$ converges in law to
$G_2(\phi_p)\sqrt{W_\infty}$
 with $G_2(\phi_p)\sim \mathcal{N}(0, A)$;
$e^{-(\alpha/2)t}\langle f_{(l)}, X_t\rangle$ converges in law to
$G_1(f_{(l)})\sqrt{W_\infty}$.

\subsubsection{Further results in the large branching rate case}

In this subsection we give two central limit theorems for the case $\alpha>2\gamma(f)b$.
Define
\begin{equation}\label{8.47}
  H_\infty:=\sum_{\gamma(f)\leq m<\alpha/(2b)}\sum_{|p|=m}a_pH^p_\infty.
\end{equation}
Let
\begin{eqnarray}
  \beta_{f_{(s)}}^2:= A\sum_{\gamma(f)\leq m<\alpha/(2b)}\frac{1}{\alpha-2mb}\sum_{|p|=m}a_p^2,\label{sigamk}
\end{eqnarray}
In Section 3.3 we will see that $\beta_{f_{(s)}}^2=\langle{\bf V}ar_{\delta_x}H_\infty,~\varphi\rangle$.

\begin{thrm}\label{The:2.1}
If $f\in \mathcal{P}$ satisfies $\alpha > 2\gamma(f)b$ and $f_{(c)}=0$, then $\sigma_{f_{(l)}}^2<\infty$.
Under $\P_{\mu}(\cdot\mid \mathcal{E}^c)$, it holds that, as $t\to\infty$,
\begin{equation}\label{the2.1}
  \left(e^{-\alpha t}\|X_t\|,~\|X_t\|^{-1/2}\left(\langle f,X_t\rangle-\sum_{\gamma(f)=m<\alpha/(2b)}e^{(\alpha-mb)t}\sum_{|p|=m}a_pH^p_\infty\right) \right)\stackrel{d}{\rightarrow}(W^*,~G_3(f)),
\end{equation}
where $W^*$ has the same distribution as $W_\infty$ conditioned on $\mathcal{E}^c$,
and $G_3(f)\sim \mathcal{N}(0,\sigma_{f_{(l)}}^2+\beta_{f_{(s)}}^2)$.
Moreover, $W^*$ and $G_3(f)$ are independent.
\end{thrm}

\begin{remark}
If $\alpha>2|p|b$, then under $\P_{\mu}(\cdot\mid \mathcal{E}^c)$, it holds that, as $t\to\infty$,
\begin{equation}\label{the:2.9}
  \left(e^{-\alpha t}\|X_t\|,~\frac{\left(\langle \phi_p,X_t\rangle-e^{(\alpha-|p|b)t}H^p_\infty\right)}{\|X_t\|^{1/2}} \right)\stackrel{d}{\rightarrow}(W^*,~G_3),
\end{equation}
where $G_3\sim\mathcal{N}(0,\frac{A}{\alpha-2|p|b})$. In particular, for $|p|=0$, we have
\begin{equation*}
  \left(e^{-\alpha t}\|X_t\|,~\frac{\|X_t\|-e^{\alpha t}W_\infty}{\sqrt{\|X_t\|}}\right)\stackrel{d}{\rightarrow}(W^*,~G_3), \quad t\to\infty,
\end{equation*}
where $G_3\sim\mathcal{N}(0,\frac{A}{\alpha})$
\end{remark}

\begin{remark} Using the theorem above, we get that if $\alpha > 2\gamma(f)b$ and $f_{(c)}=0$, then
under $\mathbb{P}_{\mu}$, we have, as $t\to\infty$,
\begin{eqnarray*}
  &&\left(e^{-\alpha t}\|X_t\|,~e^{-(\alpha/2) t}\left(\langle
 f,X_t\rangle-\sum_{\gamma(f)\le m<\alpha/(2b)}e^{(\alpha-mb)t}\sum_{|p|=m}a_pH^p_\infty\right)
  \right)\stackrel{d}{\rightarrow}(W_\infty,~\sqrt{W_\infty}G_3(f)),
\end{eqnarray*}
where $G_3(f)$ is the same as in the theorem above.
\end{remark}

\begin{thrm}\label{The:2.3}
If $f\in \mathcal{P}$ satisfies $f_{(c)}\ne 0$, then, under $\P_{\mu}(\cdot\mid \mathcal{E}^c)$,
it holds that, as $t\to\infty$,
\begin{equation}\label{the:2.3}
  \left(e^{-\alpha t}\|X_t\|,~t^{-1/2}\|X_t\|^{-1/2}\left(\langle
 f,X_t\rangle-\sum_{\gamma(f)\le m<\alpha/(2b)}e^{(\alpha-mb)t}\sum_{|p|=m}a_pH^p_\infty\right)
   \right)\stackrel{d}{\rightarrow}(W^*,~G_4(f)),
\end{equation}
where $W^*$ has the same distribution as $W_\infty$ conditioned on $\mathcal{E}^c$,
and $G_4(f)\sim \mathcal{N}(0,A\sum_{|p|=\alpha/2b}(a_p)^2)$.
Moreover, $W^*$ and $G_4(f)$ are independent.
\end{thrm}

\begin{remark}
Note that the limiting normal random variables in our Theorems  \ref{The:2.1} and \ref{The:2.3}
are non-degenerate.
\end{remark}

\section{Preliminary}

\subsection{Backbone decomposition of super-OU processes}

In this subsection,  we recall the backbone decomposition of \cite{BKS}.
Define another branching mechanism $\psi^*$ by
\begin{eqnarray}\label{psistar}
                            \psi^*(\lambda)&=& \psi(\lambda+\lambda^*)\nonumber\\
                            &=& \alpha^*\lambda+\beta\lambda^2+\int_{(0, \infty)}(e^{-\lambda x}-1+\lambda x)e^{-\lambda^* x}n(dx),
\end{eqnarray}
where
\begin{equation*}
\alpha^*=-\alpha+2\beta\lambda^*+\int_{(0, \infty)}(1-e^{-\lambda^* x})xn(dx).
\end{equation*}
It is easy to see that $\alpha^*=(\psi^*)'(0+)=\psi'(\lambda^*)> 0$.
So the $(\xi, \psi^*)$-superprocess is subcritical.
Note that it follow from \eqref{psistar} that the measure $n^*$ associated with $\psi^*$ is
$e^{-\lambda^* x}n(dx)$, thus for any $n\in\N$, $\int^\infty_0x^nn^*(dx)<\infty$.
It follows from \cite[Lemma 2]{BKS} that the
$(\xi, \psi)$-superprocess conditioned on $\mathcal{E}$ has the same law as the $(\xi, \psi^*)$-superprocess.
Let $\P^*_\mu$ be the law of the $(\xi,\psi^*)$-superprocess with initial configuration $\mu$, and define
\begin{equation*}
  u^*_f(x,t)=-\log \P^*_{\delta_x}(e^{-\langle f,X_t\rangle}).
\end{equation*}

It was shown in \cite{E.B2.} that one can associate with $\{\P^*_{\delta_x}:x\in \mathbb{R}^d\}$ a family
of measures $\{\N^*_x:x\in \mathbb{R}^d\}$, defined on the same measurable space as the probabilities
$\{\mathbb{P}^*_{\delta_x}:x\in \R^d\}$
and satisfying
\begin{equation}
\mathbb{N}^{*}_x (1- e^{-\langle f, X_t \rangle}) = -\log \mathbb{P}^{*}_{\delta_x}(e^{-\langle  f, X_t\rangle}) = u^{*}_f(x,t),
\label{DK}
\end{equation}
for all $f\in {\cal B}^+_b(\R^d)$
and $t\geq 0$. Intuitively speaking, the branching property implies that $\mathbb{P}^*_{\delta_x}$
is an infinitely divisible measure on the path space of $X$, that is to say, the space of measure-valued
cadlag functions, $ \mathbb{D}([0,\infty)\times \mathcal{M}_F({\R^d}))$, and  (\ref{DK})  is a
`L\'evy-Khinchine' formula in which  $\mathbb{N}^*_x$ plays the role of its `L\'evy measure'.
The measures $\{\N^*_x:x\in \mathbb{R}^d\}$ will play a crucial role in the forthcoming analysis.

Let $\mathcal{M}_a(\mathbb{R}^d)$ be the space of finite atomic measures on $\mathbb{R}^d$. For $\nu\in\mathcal{M}_a(\R^d)$,
let $Z=(Z_t:t\geq 0)$ be a branching OU-process with initial configuration $\nu$. $\{Z_t,t\geq 0\}$ is an $\mathcal{M}_a(\mathbb{R}^d)$-valued Markov process in which individuals, from the moment of birth,
live for an independent and exponential distributed period of time with parameter $\alpha^*$ during
which they move according to the OU-process issued from their position of birth and
at death they give birth at the same position to an independent number of
offspring with distribution  $(p_n:n\geq0)$, where $p_0=p_1=0$ and for $n\geq 2$,
\begin{equation*}
  p_n=\frac{1}{\lambda^*\alpha^*}\left\{\beta(\lambda^*)^2
  \mathbf{1}_{\{n=2\}}+(\lambda^*)^n\int_{(0,\infty)}\frac{x^n}{n!}e^{-\lambda^*x}n(dx)\right\}.
\end{equation*}
The generator of $Z$ is given by
\begin{eqnarray}\label{2.16}
  F(s)&=&\alpha^* \sum_{n\geq0}p_n(s^n-s)=\frac{1}{\lambda^*}\psi(\lambda^*(1-s)).
\end{eqnarray}
$Z$ is refereed as the $(\xi, F)$-backbone in \cite{BKS}.
Moreover, when referring to individuals in $Z$ we will use the classical Ulam-Harris notation so that every particle in
$Z$ has a unique label, see \cite{HH}.
Let ${\cal T}$ be the set of labels of individuals realized in $Z$.
Let  $|Z_t|$ be the number of particles alive at time $t$.
 For each individual
$u\in{\cal T}$ we shall write
$\tau_u$ and $\sigma_u$ for its birth and death times respectively and
$\{z_u(r): r\in [\tau_u,\sigma_u]\}$ for its spatial trajectory.
Now we describe three kinds of immigrations along the backbone $Z$ as follows.

\begin{enumerate}
  \item {\bf Continuous immigration:} The process $I^{\N^*}$ is measure-valued on $\R^d$ such that
  $$
  I^{\N^*}_t:=\sum_{u\in{\cal T}}\sum_{u\wedge\tau_u<r\le t\wedge\sigma_u}X^{(1, u, r)}_{t-r}\ ,
  $$
    where, given $Z$, independently for each $u\in{\cal T}$ with $\tau_u<t$, the processes
  $X^{(1, u, r)}_{\cdot}$ are  independent copies of the canonical process $X$,
   immigrated along the space-time trajectory $\{(z_u(r), r): r\in (\tau_u, t\wedge\sigma_u]\}$ with rate
  $2\beta{\rm d} r\times{\rm d}\N^*_{z_u(r)}$.

   \item {\bf Discontinuous immigration:} The processes $I^{\P^*}$ is measure-valued on $\R^d$ such that
  $$
  I^{\P^*}_t:=\sum_{u\in{\cal T}}\sum_{t\wedge\tau_u<r\le t\wedge\sigma_u}X^{(2,u,r)}_{t-r}\ ,
  $$
    where, given $Z$, independently for each $u\in{\cal T}$ with $\tau_u<t$,
    the processes $X^{(2, u, r)}_{\cdot}$ are independent copies of the canonical process $X$,
  immigrated along the space-time trajectory $\{z_u(r): r\in(\tau_u, t\wedge\sigma_u]\}$
  with rate ${\rm d}r\times\int_{y\in(0,\infty)}ye^{-\lambda^* y}n({\rm }dy)\times{\rm d}
   \P^*_{y\delta_{z_u(r)}}$.
  \item {\bf Branching point biased immigration:} The process $I^{\eta}$ is measure-valued on
  $\R^d$
  such that
  \begin{equation*}
  I^{\eta}_t=\sum_{u\in{\cal T}} \mathbf{1}_{\sigma_u\leq t}X^{(3,u)}_{t-\sigma_u}\ ,
  \end{equation*}
where, given $Z$, independently for each $u\in{\cal T}$ with ${\sigma_u}\le t$,
the processes $X^{(3, u)}_{\cdot}$ are independent copies of the canonical process $X$ issued
at time ${\sigma_u}$ with law $\mathbb{P}^*_{Y_u\delta_{z_u({\sigma_u})}}$ where,
given $u$ has $n\geq 2$ offspring, the independent random variable $Y_u$ has distribution
$\eta_n(z_u(r), {\rm d}y)$, where
      \begin{equation*}
        \eta_n(dy)=\frac{1}{p_n\lambda^*\alpha^*}\left\{\beta(\lambda^*)^2
        \delta_0(dy)\mathbf{1}_{\{n=2\}}+(\lambda^*)^n\frac{y^n}{n!}e^{-\lambda^*y}n(dy)\right\}.
      \end{equation*}
\end{enumerate}
Now we define another $\mathcal{M}_F(\mathbb{R}^d)$-valued process $I=\{I_t : t\geq 0\}$ by
\begin{equation*}
  I:=I^{\N^*}+I^{\P^*}+I^{\eta}\ ,
\end{equation*}
where the processes $I^{\N^*}=\{I^{\N^*}_t: t\geq0\}$, $I^{\P^*}=\{I^{\P^*}_t: t\geq0\}$ and
$I^{\eta}=\{I^{\eta}_t: t\geq0\}$,
conditioned on $Z$, are independent of each other.
We denote the law of $I$ by $\mathbb{Q_\nu}$. Recall that $\nu$ is the initial configuration of $Z$.

For $\mu \in \mathcal{M}_F(\mathbb{R}^d)$,
let $\widetilde{X}$ be an independent copy of $X$ under $\P^*_\mu$
and be independent of $I$.
Then we define a measure-valued process $\Lambda=\{\Lambda_t: t\geq0\}$ by
\begin{equation}\label{12}
  \Lambda=\widetilde{X}+I.
\end{equation}
Note that $Z$, $\widetilde{X}$ and the three immigration processes above
are defined on the same probability space.
 We denote the law of $\Lambda$ by ${\bf P}_{\mu\times\nu}$.
When $\nu$ is a Poisson random measure with intensity
measure $\lambda^*\mu$, then we write this law by
${\bf P}_\mu$.
The following result is proved in \cite{BKS}.

\begin{prop}
 For any $\mu \in\mathcal{M}_F(\mathbb{R}^d)$, the process $(\Lambda, \bP_\mu)$ is Markovian and has the same law as $(X, \P_\mu)$.
\end{prop}

We will need the following $\sigma$-fields later on:
\begin{eqnarray}
\mathcal{F}_t&=&\sigma(\Lambda_s, s\leq t), \quad t\ge 0,\label{e:calFt}\\
\mathcal{G}_t&=&\sigma(\Lambda_s, Z_s,s\leq t), \quad t\ge 0.\label{e:calGt}
\end{eqnarray}

\subsection{Moments}

Now we use Laplace transforms to calculate the moments of $X$.
We will omit some details, for these omitted details, see \cite{E.B.}.
For any $f\in {\cal P}$, we define
 \begin{equation*}
    u_f(x,t,\theta)=-\log \P_{\delta_x}(e^{-\langle\theta f,X_t\rangle}),
 \end{equation*}
then
\begin{equation}\label{2.13}
  u_f(x,t,\theta)+\Pi_x\int_0^t\psi(u_f(\xi_s,t-s,\theta))ds=\theta\Pi_x f(\xi_t).
\end{equation}
Differentiating both sides of \eqref{2.13} with respect to $\theta$, we get
\begin{eqnarray}
     u_f^{(1)}(x,t,0) &=& e^{-\psi'(0+)t}T_tf(x), \label{1.6}\\
     u_f^{(2)}(x,t,0) &=& -\psi^{(2)}(0+)\int_0^t e^{-\psi'(0+)(t-s)}T_{t-s} [u_f^{(1)}(\cdot,s,0)]^2(x)ds\nonumber\\
                      &=&-Ae^{\alpha t}\int_0^te^{\alpha s}T_{t-s}[T_sf]^2(x)\,ds.\label{2.14}
\end{eqnarray}
Let $\mu_t:=\P_{\mu}\langle f,X_t\rangle$. The moments are given by
\begin{equation*}
  \P_{\mu}(\langle f,X_t\rangle)^{n}=(-1)^n(e^{-\langle u_f,\mu\rangle})^{(n)}|_{\theta=0}.
\end{equation*}
In particular,
\begin{equation}
  \mu_t=\P_{\mu}\langle f,X_t\rangle = \langle u_f^{(1)}(x,t,0),\mu\rangle=e^{\alpha t}\langle T_tf,\mu\rangle,\label{2.17}\end{equation}
\begin{equation}
  \P_{\mu}(\langle f,X_t\rangle-\mu_t)^2 =-\langle u_f^{(2)}(x,t,0),\mu\rangle.\label{2.15}
\end{equation}

Recall that $\widetilde{X}_t$ is defined in Section 2.1.
It is a subcritical superprocess with branching mechanism $\psi^*(\lambda)=\psi(\lambda+\lambda^*)$. Thus $(\psi^*)^{(m)}(0+)=\psi^{(m)}(\lambda^*)$ exists for all $m\geq 1$.
For any $f\in {\cal P}$, we define
 \begin{equation}\label{4.14}
     u_f^*(x,t,\theta)=-\log \P_{\delta_x}(e^{-\langle\theta f,\widetilde{X}_t\rangle}).
 \end{equation}
Then
\begin{equation}\label{4.13}
  u^*_f(x,t,\theta)+\Pi_x\int_0^t\psi^*(u_f^*(\xi_s,t-s,\theta))ds=\theta\Pi_x f(\xi_t).
\end{equation}
Differentiating both sides of \eqref{4.14} with respect to $\theta$, we have
\begin{eqnarray}
      (u_f^*)^{(1)}(x,t,0) &=& e^{-\alpha^*t}T_tf(x), \label{4.17}\\
     (u^*_f)^{(2)}(x,t,0) &=& -(\psi^*)^{(2)}(0+)\int_0^t e^{-(\psi^*)'(0+)(t-s)}T_{t-s} [(u^*_f)^{(1)}(\cdot,s,0)]^2(x)ds\nonumber\\
                      &=&-(\psi^*)^{(2)}(0+)e^{-\alpha^* t}\int_0^te^{-\alpha^* s}T_{t-s}[T_sf]^2(x)\,ds,\label{24.14}\\
     (u_f^*)^{(3)}(x,t,0) &=& -(\psi^*)^{(3)}(0+)\int_0^te^{-\alpha^* s}T_{s}[ (u_f^*)^{(1)}(\cdot,t-s,0)]^3(x)\,ds\nonumber\\
     &&-3(\psi^*)^{(2)}(0+)\int_0^t e^{-\alpha^* s}T_s[( (u_f^*)^{(1)}(u_f^*)^{(2)})(\cdot,t-s,0)](x)\,ds, \label{4.18}
\end{eqnarray}
and
\begin{eqnarray}
     (u_f^*)^{(4)}(x,t,0) = -\int_0^te^{-\alpha^* s}T_{s}[ J(\cdot,t-s)](x)\,ds, \phantom{3in}\label{4.19}
   \end{eqnarray}
where
\begin{eqnarray*}
  J(x,t)&=& \left[(\psi^*)^{(4)}(0)\left((u_f^*)^{(1)}\right)^4+6(\psi^*)^{(3)}(0)\left((u_f^*)^{(1)}\right)^2(u_f^*)^{(2)}\right](x,t,0)\nonumber\\
  &&+\left[4(\psi^*)^{(2)}(0)(u_f^*)^{(1)}(u_f^*)^{(3)}+3(\psi^*)^{(2)}(0)\left((u_f^*)^{(2)}\right)^2\right](x,t,0).
\end{eqnarray*}
By \eqref{4.14}, the moments of $\widetilde X$ are given by
\begin{equation*}
  \P_{\mu}(\langle f,\widetilde{X}_t\rangle)^{n}=(-1)^n(e^{-\langle u_f^*,\mu\rangle})^{(n)}|_{\theta=0}.
\end{equation*}
In particular, we have
\begin{equation}
  \P_{\mu}\langle f,\widetilde X_t\rangle = \langle (u_f^*)^{(1)}(x,t,0),\mu\rangle=e^{-\alpha^* t}\langle T_tf,\mu\rangle,\label{2.17*}\end{equation}
\begin{equation}
  \P_{\mu}(\langle f,\widetilde X_t\rangle-\P_{\mu}\langle f,\widetilde X_t\rangle)^2 =-\langle (u^*_f)^{(2)}(x,t,0),\mu\rangle.\label{2.15*}
\end{equation}
\begin{eqnarray}
   \P_{\mu}(\langle f,\widetilde{X}_t\rangle-\P_{\mu}\langle f,\widetilde{X}_t\rangle)^4= -\langle (u_f^*)^{(4)}(x,t,0),\mu\rangle+
   3\langle (u_f^*)^{(2)}(x,t,0),\mu\rangle^2\label{4.15*}.
\end{eqnarray}

\subsection{Estimates on the semigroup $T_t$}

Recall that $\xi=\{\xi_t: t\ge 0\}$ is the OU process and $\{T_t\}$ is the semigroup of $\xi$. It is well known that
under $\Pi_x$, $\xi_t\sim\mathcal{N}(xe^{-b t}, \sigma_t^2)$,
where $\sigma_t^2=\sigma^2(1-e^{-2b t})/(2b)$.
Let $G$ be an $\mathbb{R}^d$-valued standard normal  random variable,
then using $(a+b)^n\leq 2^n(a^n+b^n),~a\geq0,~b\geq0$, we get
\begin{equation}\label{2.35}
  T_t(\|\cdot\|^n)(x)=E(\|\sigma_tG+xe^{-bt}\|^n)\leq 2^n\left[(\sigma/\sqrt{2b})^nE(\|G\|^n)+\|x\|^n\right].
\end{equation}
Using this, we can easily get that
\begin{equation}\label{e:ousmgp}
T_t(1+\|\cdot\|^n)(x)\le c(n)(1+\|x\|^n),
\end{equation}
where $c(n)$ does not depend on $t$.
\begin{lemma}\label{lem:2.1}
For any $f\in L^2(\varphi)$, we have that, for every $x\in \mathbb{R}^d$,
\begin{equation}\label{2.8}
  T_tf(x)=\sum_{n=\gamma(f)}^{\infty}e^{-nb t}\sum_{|p|=n}a_p\phi_p(x),
\end{equation}
\begin{equation}\label{2.6}
  \lim _{t\to \infty}e^{\gamma(f) bt}T_tf(x)=\sum_{|p|=\gamma(f)}a_p\phi_p(x).
\end{equation}
Moreover, there exists $c>0$ such that for $t\geq 1$,
\begin{equation}\label{2.7}
|T_tf(x)|\leq c e^{-\gamma(f) bt}e^{\frac{b}{2\sigma^2}\|x\|^2}, \quad x\in \mathbb{R}^d.
\end{equation}
\end{lemma}
\textbf{Proof:}\quad For every $f\in L^2(\varphi)$, using the fact that
$\varphi(x)$ is the invariant density of $\xi$ we get that
\begin{equation}\label{3.8}
  \int \varphi(x)\left(T_t|f|(x)\right)^2\,dx\leq \int \varphi(x)T_t[|f|^2](x)\,dx=\int |f(y)|^2\varphi(y)\,dy<\infty,
\end{equation}
 so $T_tf(x)\in L^2(\varphi)$. Moreover, by the fact $\xi_t\sim\mathcal{N}(xe^{-b t}, \sigma_t^2)$,
 $T_t|f|(x)$ is continuous in $x$. Thus $T_t|f|(x)<\infty$ for all $x\in\mathbb{R}^d$. \eqref{3.8}
implies that $T_t$ is a bounded linear operator on $L^2(\varphi)$.
Let $f_k(x)=\sum_{n=0}^{k}\sum_{|p|=n}a_p\phi_p(x)$. Since $f_k\to f$ in $L^2(\varphi)$, we have
$T_tf_k\to T_tf$ in $L^2(\varphi)$, as $k\to\infty$. By linearity, we have
\begin{equation*}
  T_tf_k(x)=\sum_{n=0}^{k}e^{-nb t}\left(\sum_{|p|=n}a_p\phi_p(x)\right).
\end{equation*}
We claim that the series $\sum_{n=0}^{\infty}e^{-nbt}\left(\sum_{|p|=n}a_p\phi_p(x)\right)$ is uniformly convergent
on any compact subset of $\mathbb{R}^d$. Thus
$\sum_{n=0}^{\infty}e^{-nbt} \left(\sum_{|p|=n}a_p\phi_p(x)\right)$ is continuous in $x$. So for all $x\in\mathbb{R}^d$,
\begin{equation*}
  T_tf(x)=\sum_{n=0}^{\infty}e^{-nbt}\left(\sum_{|p|=n}a_p\phi_p(x)\right).
\end{equation*}
Now we prove the claim. In fact, by Cramer's inequality (for example, see \cite[Equation (19) on p.\/207]{EMOT}),
for all $p\in \mathbb{Z}^d_+$ we have
\begin{equation}\label{e:cramer}
   |\phi_p(x)|\leq K e^{\frac{b}{2\sigma^2}\|x\|^2},
\end{equation}
where $K$ is a constant. So we only need to prove
$\sum_{n=0}^{\infty}e^{-nbt}\left(\sum_{|p|=n}|a_p|\right)<\infty$.
By H\"older's inequality,
\begin{equation}\label{2.11}
\sum_{n=0}^{\infty}e^{-nbt}\left(\sum_{|p|=n}|a_p|\right) \leq\left(\sum_{n=\gamma(f)}^{\infty}K_ne^{-2nb t}\right)^{1/2}\left(\sum_{n=\gamma(f)}^{\infty}\sum_{|p|=n}|a_p|^2\right)^{1/2},
\end{equation}
where $K_n={n+d-1 \choose d-1}=\sharp\{p\in\mathbb{Z}_+^d: |p|=n\}$.
Since $K_n\leq (n+d)^d$, we have that $\sum_{n=\gamma(f)}^{\infty}K_ne^{-2nbt}<\infty$.
Using the fact that $\{\phi_p(x), ~p\in \mathbb{Z}_+^d\}$ form a complete orthogonal basis for $L^2(\varphi)$,
we get $\sum_{n=\gamma(f)}^{\infty}\sum_{|p|=n}|a_p|^2=\int \varphi(x)|f(x)|^2\,dx<\infty$. Therefore the claim is
true.

By \eqref{e:cramer} and \eqref{2.11}, for $t\geq 1$, we have
\begin{eqnarray}
  |T_tf(x)|&\leq & e^{-\gamma(f) bt}\left(\sum_{n=0}^{\infty}K_{n+\gamma(f)}e^{-2nb}\right)^{1/2}
  (\sum_{n=\gamma(f)}^{\infty}\sum_{|p|=n}|a_p|^2)^{1/2}Ke^{\frac{b}{2\sigma^2}\|x\|^2}
  \label{3.18}\\
   &\lesssim & e^{-\gamma(f) bt}e^{\frac{b}{2\sigma^2}\|x\|^2}, \quad x\in \R^d.
  \nonumber
\end{eqnarray}
Therefore, for $t\geq 1$,
\begin{eqnarray}
|e^{\gamma(f) bt}T_tf(x)-\sum_{|p|=\gamma(f)}a_p\phi_p(x)|&=& e^{\gamma(f) bt}|T_tf(x)-e^{-\gamma(f) bt}\sum_{|p|=\gamma(f)}a_p\phi_p(x)|\nonumber\\
&=& e^{\gamma(f) bt}\left|T_t(f-\sum_{|p|=\gamma(f)}a_p\phi_p)(x)\right|\nonumber\\
&\lesssim& e^{-bt}e^{\frac{b}{2\sigma^2}\|x\|^2},
\end{eqnarray}
which implies \eqref{2.6}. The proof is now complete. \hfill$\Box$

\smallskip

For $p\in \mathbb{Z}_+^d$, we use the notation
$f^{(p)}(x):=\frac{\partial}{\partial x_1^{p_1}\partial x_2^{p_2}\cdots\partial x_d^{p_d}}f(x)$. Define
\begin{equation*}
  \mathcal{P}^*=\{f\in C^\infty: ~f^{(p)}\in \mathcal{P} \mbox{ for all }p\in \mathbb{Z}_+^d\}.
\end{equation*}
It can be easily shown that, for any $f\in\mathcal{P}$,  $T_tf(x)\in \mathcal{P}^*$.

\begin{lemma}
For any $f\in \mathcal{P}^*$ and $p\in\mathbb{Z}_+^d$
satisfying $0\leq |p|\leq \gamma(f)$, we have $\gamma(f^{(p)})\geq\gamma(f)-|p|$.
\end{lemma}
\textbf{Proof:}\quad By the definition of $\phi_p$ and $\varphi$, it is easy to check that
\begin{equation*}
  \phi_p(x)\varphi(x)=(-1)^{|p|}c_p\varphi^{(p)}(x),
\end{equation*}
where $c_p=\frac{1}{\sqrt{p!2^{|p|}}}\left(\frac{\sigma^2}{b}\right)^{|p|/2}$. Integrating by parts, we get
\begin{equation}\label{1.4}
  \int f(x)\phi_p(x)\varphi(x)\,dx=c_p\int_{\mathbb{R}^d} f^{(p)}(x)\varphi(x)\,dx.
\end{equation}
Thus
\begin{equation*}
  \gamma(f^{(p)})=
  \inf\{k: \mbox{ there~exists } p \mbox{ such that } |p|=k \mbox{ and }
  \int_{\mathbb{R}^d} f^{(p)}(x)\varphi(x)\,dx\ne 0\}.
\end{equation*}
Hence if $|p'|<\gamma(f)-|p|$, we have $\int_{\mathbb{R}^d} f^{(p+p')}(x)\varphi(x)\,dx= 0$,
which implies $\gamma(f^{(p)})\geq \gamma(f)-|p|$. \hfill$\Box$

\smallskip

In the following lemma, we give another estimate for $T_tf$, which will be very useful later.
\begin{lemma}\label{lem:2.3}
For every $f\in \mathcal{P}$, there exist $r\in\N$ and $c>0$ such that
\begin{eqnarray}
  &&e^{\gamma(f) bt}|T_tf(x)|\leq c (1+\|x\|^{r}),\label{2.40}\\
  &&\left|e^{\gamma(f) bt}T_tf(x)-\sum_{|p|=\gamma(f)}a_p\phi_p(x)\right|\leq ce^{-bt} (1+\|x\|^{r})\label{2.41}.
\end{eqnarray}
\end{lemma}
\textbf{Proof:}\quad
Let $g(x)=T_1f(x)$$\in \mathcal{P}^*$. Then $\gamma(g)=\gamma(f)$ and
there exist $k\in \N$ and $c_1>0$ such that,
for $|p|=0,1\cdots,\gamma(f) $,
$|g^{(p)}(x)|\leq c_1(1+\|x\|^k)$. For $x=(x_1,x_2\cdots,x_d)\in \mathbb{R}^d$,
we define $x^p:=\prod_{i=1}^d x_i^{p_i}$.
Then for $s>0$ we have
\begin{eqnarray*}
  T_sg(x)&=&T_s[g(\cdot+xe^{-bs})](0)\\
  &=&T_s\left[g(\cdot+xe^{-bs})-\sum_{m=0}^{\gamma(f)-1}\sum_{|p|=m}g^{(p)}(\cdot)
  x^{p}e^{-mbs}/p!\right](0)\\
  &&+\sum_{m=0}^{\gamma(f)-1}\sum_{|p|=m}T_s[g^{(p)}](0)x^{p}e^{-mbs}/p!\\
  &=& (I)+(II).
\end{eqnarray*}
It follows from \eqref{2.7} and the fact that $\gamma(g^{(p)})\geq\gamma(g)-|p|$,
we have 
$$\sup_{s>0}e^{(\gamma(g)-|p|)b s}|T_s[g^{(p)}](0)|<\infty.
$$ 
Thus
 \begin{equation*}
   |(II)|\lesssim e^{-\gamma(f) bs}\sum_{m=0}^{\gamma(f)-1}\sum_{|p|=m}|x^p|\lesssim e^{-\gamma(f) bs}(1+\|x\|^{\gamma(f)}).
 \end{equation*}
 Using Taylor's formula and the fact $|g^{(p)}(x)|\lesssim 1+\|x\|^k$, we get
 \begin{eqnarray*}
  \left| g(y+xe^{-bs})-\sum_{m=0}^{\gamma(f)-1}\sum_{|p|=m}g^{(p)}(y)x^{p}e^{-mbs}/p!\right|
  &=&\sum_{|p|= \gamma(f)}|g^{(p)}(\theta)||x^p| e^{-\gamma(f) bs}/(\gamma(f) !)\\
  &\lesssim& (1+\|y\|^k+\|x\|^k)
  |x|^{\gamma(f)}e^{-\gamma(f) bs},
 \end{eqnarray*}
 where $\theta$ is a point on the line segment connecting $y$ and $y+xe^{-bs}$.
 Then by the fact that $T_s[\|\cdot\|^k](x)\lesssim
   1+\|x\|^k$, we get $\sup_{s>0}T_s[\|\cdot\|^k](0)<\infty$.
 Therefore, we have
 \begin{equation*}
   |(I)|\lesssim (1+\|x\|^{k+\gamma(f)})e^{-\gamma(f) bs}.
 \end{equation*}
 Consequently,
 \begin{equation*}
  e^{\gamma(f) bs}|T_s g|(x)\lesssim 1+\|x\|^{k+\gamma(f)}.
 \end{equation*}
Let $r_1=k+\gamma(f)$. For $t\geq 1$, combining $T_tf(x)=T_{t-1}(g)(x)$
with the above inequality,
we arrive at \eqref{2.40} for $t\ge 1$. For $t<1$,
\begin{equation*}
  e^{\gamma(f) bt}|T_tf(x)|\lesssim e^{\gamma(f) b}(1+\|x\|^k)\lesssim 1+\|x\|^{r_1},
\end{equation*}
so \eqref{2.40} is also valid.

It follows from \eqref{2.40} that there exists $r_2\in\N$ such that
\begin{equation*}
  e^{(\gamma(f)+1)bt}\left|T_tf(x)-e^{-\gamma(f) bt}\sum_{|p|=\gamma(f)}a_p\phi_p(x)\right|
  \lesssim 1+\|x\|^{r_2}.
\end{equation*}
Now \eqref{2.41} follows immediately. \hfill$\Box$

\smallskip

From the above calculations, we have
\begin{lemma} \label{lem:2.2}
Let $f\in\mathcal{P}$.
\begin{description}
\item{(i)} If $\alpha<2\gamma(f) b$, then
\begin{equation*}
\lim_{t\to \infty}e^{-(\alpha/2)t}\P_{\delta_x}\left(\langle f, X_t\rangle\right)=0,
\end{equation*}
\begin{equation}\label{2.2}
\lim_{t\to \infty}e^{-\alpha t}\mathbb{V}ar_{\delta_x}\langle f,X_t\rangle = \sigma_f^2,
\end{equation}
where $\mathbb{V}ar_{\delta_x}$ stands for the variance under $\P_{\delta_x}$
and $\sigma_f^2$ is defined in \eqref{e:sigmaf}.
\item{(ii)} If $\alpha=2\gamma(f) b$, then
\begin{equation}\label{4.10}
\lim_{t\to \infty}t^{-1/2} e^{-(\alpha/2)t}\P_{\delta_x}\left(\langle f, X_t\rangle\right)=0,
 \end{equation}
and there exists $r\in \N$ such that
\begin{equation}\label{1.18}
  |t^{-1}e^{-\alpha t}\mathbb{V}ar_{\delta_x}\langle f,X_t\rangle|\lesssim 1+\|x\|^{2r}
\end{equation}
and
\begin{equation}\label{1.17}
|t^{-1}e^{-\alpha t}\mathbb{V}ar_{\delta_x}\langle f,X_t\rangle-\rho_f^2|\lesssim t^{-1}(1+\|x\|^r),
\end{equation}
which in particular implies that
\begin{equation}\label{3.11}
\lim_{t\to \infty}t^{-1}e^{-\alpha t}\mathbb{V}ar_{\delta_x}\langle f,X_t\rangle = \rho_f^2,
\end{equation}
where $\rho^2_f$ is defined in \eqref{1.3}.
\item{(iii)} If $\alpha > 2\gamma(f) b$, then
\begin{equation}\label{3.12}
\lim_{t\to \infty}e^{-2(\alpha-\gamma(f) b)t}\mathbb{V}ar_{\delta_x}\langle f,X_t\rangle=\eta_f^2(x),
\end{equation}
where
\begin{equation}\label{e:etaf}
\eta_f^2(x)=A\int_0^\infty e^{-(\alpha-2\gamma(f) b)s}T_s\left(\sum_{|p|=\gamma(f)}a_p\phi_p\right)^2(x)\,ds.
\end{equation}
\end{description}
\end{lemma}
\textbf{Proof:}\quad
It follows from \eqref{2.14} and \eqref{2.15} that
\begin{equation}\label{2.24}
\mathbb{V}ar_{\delta_x}\langle f,X_t\rangle=Ae^{\alpha t}\int_0^t e^{\alpha s}T_{t-s} [T_sf]^2(x)ds
=Ae^{2\alpha t}\int_0^t e^{-\alpha s}T_{s} [T_{t-s}f]^2(x)ds.
\end{equation}

(i)
If $\alpha<2\gamma(f) b$, by Lemma~\ref{lem:2.1},  we have $\lim _{t\to \infty}e^{\gamma(f) bt}T_tf(x)=\sum_{|p|=\gamma(f)}a_p\phi_p(x)$. Thus
\begin{equation*}
\lim_{t\to\infty}e^{-(\alpha/2) t}\P_{\delta_x}\langle f,X_t\rangle=\lim_{t\to\infty}e^{(\alpha-2\gamma(f) b)t/2}[e^{\gamma(f) b t}T_tf(x)]= 0.
\end{equation*}
It follows from Lemma~\ref{lem:2.3} that there exists  $r\in \N$ such that
$e^{\gamma(f) bs}|T_sf|(x)\lesssim 1+\|x\|^r$. Using \eqref{e:ousmgp}, we have
\begin{equation}\label{3.10}
  T_{t-s}[e^{\gamma(f) bs}T_sf]^2(x)\lesssim 1+\|x\|^{2r}.
\end{equation}
Thus $e^{\alpha s}T_{t-s} [T_sf]^2(x)\lesssim e^{(\alpha-2\gamma(f) b)s}(1+\|x\|^{2r})$.
Hence by the dominated convergence theorem, we get
\begin{equation*}
\lim_{t\to\infty}\int_0^te^{\alpha s}T_{t-s}[T_sf]^2(x)\,ds=
 \int_0^{\infty} e^{\alpha s}\langle(T_sf)^2,\varphi\rangle\, ds,
\end{equation*}
which implies \eqref{2.2}.

(ii)
If $\alpha=2\gamma(f) b$, then by \eqref{2.24}, we have
\begin{equation}\label{2.43}
t^{-1}e^{-\alpha t}\mathbb{V}ar_{\delta_x}\langle f,X_t\rangle=At^{-1}\int_0^tT_{t-s}[e^{\gamma(f)bs}T_sf]^2(x)\,ds.
\end{equation}
By Lemma \ref{lem:2.3}, there exists $r\in \N$ satisfying
\eqref{2.40}, \eqref{2.41} and
\begin{equation*}
 |\sum_{|p|=\gamma(f)}a_{p}\phi_{p}(x)|\lesssim 1+\|x\|^{r},
\end{equation*}
which follows from the fact that $\phi_p(x)$ is a polynomial. Then by \eqref{2.40} and \eqref{2.43}, it is easy to get \eqref{1.18}.

Let $h(x):=(\sum_{|p|=\gamma(f)}a_{p}\phi_{p}(x))^2$.
Then we have
\begin{eqnarray*}
  &&|(e^{\gamma(f) bs}T_sf(x))^2-h(x)|\\
  &\leq&\left|e^{\gamma(f) bs}T_sf(x)-\sum_{|p|=\gamma(f)}a_{p}\phi_{p}(x)\right|
  \left(e^{\gamma(f) bs}|T_sf|(x)+|\sum_{|p|=\gamma(f)}a_{p}\phi_{p}(x)|\right)\\
  &\lesssim& e^{-bs}(1+\|x\|^{2r}).
\end{eqnarray*}
Since $\gamma(h)=0$ and $\sum_{|p|=\gamma(f)}a_p^2=\langle h, \varphi\rangle$, by \eqref{2.41},
there exists $r'\in \N$ such that
\begin{equation}\label{5.7}
  \left|T_{t-s}h(x)-\sum_{|p|=\gamma(f)}a_p^2\right|\lesssim e^{-b(t-s)}(1+\|x\|^{r'}).
\end{equation}
Let $r_0=\max(2r,r')$, then
\begin{eqnarray*}
 \left| T_{t-s}(e^{\gamma(f) bs}T_sf)^2(x)-\sum_{|p|=\gamma(f)}a_p^2\right|&\leq& T_{t-s} |(e^{\gamma(f) bs}T_sf(x))^2-h(x)|+ \left|T_{t-s}h(x)-\sum_{|p|=\gamma(f)}a_p^2\right|\\
 &\lesssim &(e^{-bs}+e^{-b(t-s)})(1+\|x\|^{r_0}).
\end{eqnarray*}
It follows that
\begin{eqnarray}
&&\frac{1}{t}\int_{0}^t|T_{t-s}(e^{\gamma(f) bs}T_sf)^2(x)-\sum_{|p|=\gamma(f)}a_{p}^2|\,ds\nonumber\\
  &\lesssim &\frac{\int_{0}^{t}(e^{-bs}+e^{-b(t-s)})(1+\|x\|^{r_0})\,ds}{t}\lesssim t^{-1}(1+\|x\|^{r_0}).\label{5.10}
\end{eqnarray}
Then \eqref{1.17} follows from \eqref{2.43} and \eqref{5.10}.

(iii)
If $\alpha > 2\gamma(f) b$, then by \eqref{2.24}, we have
\begin{eqnarray*}
e^{-2(\alpha-\gamma(f) b)t} \mathbb{V}ar_{\delta_x}\langle f,X_t\rangle =
A\int_0^te^{-(\alpha-2\gamma(f) b)s}T_s[e^{\gamma(f) b(t-s)}T_{t-s}f]^2(x)\,ds.
\end{eqnarray*}
By Lemma~\ref{lem:2.3}, there exists $r\in \N$ such that
$[e^{\gamma(f) b(t-s)}T_{t-s}f(x)]^2\leq c(1+\|x\|^{2r})$.
Thus $$T_s[e^{\gamma(f) b(t-s)}T_{t-s}f]^2(x)\lesssim 1+\|x\|^{2r}.$$
Now by the dominated convergence theorem and \eqref{2.6}, we have
\begin{equation*}
 \lim_{t\to\infty} \int_0^te^{-(\alpha-2\gamma(f) b)s}T_s[e^{\gamma(f) b(t-s)}T_{t-s}f]^2(x)\,ds=
 A\int_0^\infty e^{-(\alpha-2\gamma(f) b)s}T_s\left(\sum_{|p|=\gamma(f)}a_p\phi_p\right)^2(x)\,ds.
\end{equation*}
The proof of (iii) is now complete. \hfill$\Box$

\smallskip

According to~\cite{BKS}, under ${\bf P}_{\delta_x}$,  we have that, conditioned on
$\mathcal{F}_t$ (see \eqref{e:calFt}), the backbone $Z_t$ is a Poisson point process with the intensity $\lambda^*\Lambda_t$. In particular, $Z_0=N\delta_x$, where $N$ is a Poisson random variable with
parameter $\lambda^*$. Then we have
\begin{equation}\label{2.25}
  \Lambda_t=\widetilde{X}_t+\sum_{j=1}^{N}I^{j}_t,
\end{equation}
 where $I^{j},j=1,2,...$ are independent copies of $I$ under $\mathbb{Q}_{\delta_x}$ and
 are independent of $N$. The first moment of $I$ can be calculated by
 \begin{equation}\label{e:newexp}
   {\bf P}_{\delta_x}\langle f,\Lambda_t\rangle={\bf P}_{\delta_x}\langle f,\widetilde{X}_t\rangle+\lambda^*\mathbb{Q}_{\delta_{x}}\langle f,I_t\rangle.
 \end{equation}
 Thus
 \begin{equation}\label{2.27}
   \nu_t:= \mathbb{Q}_{\delta_{x}}\langle f,I_t\rangle = \frac{1}{\lambda^*}\left({\bf P}_{\delta_x}\langle f,\Lambda_t\rangle-{\bf P}_{\delta_x}\langle f,\widetilde{X}_t\rangle\right) = \frac{1}{\lambda^*}(e^{\alpha t}-e^{-\alpha^* t})T_tf(x).
 \end{equation}
 For the second moment, let ${\bf V}ar_{\delta_x}$ stand for the variance under ${\bf P}_{\delta_x}$
  and $\mathbb{V}_{\delta_x}$ stand for the variance under $\mathbb{Q}_{\delta_x}$.
 By \eqref{2.25}, we have
 \begin{equation*}
  {\bf V}ar_{\delta_x}\langle f,\Lambda_t\rangle ={\bf V}ar_{\delta_x}\langle f,\widetilde{X}_t\rangle+\lambda^*\mathbb{Q}_{\delta_x}\langle f,I_t\rangle^2.
 \end{equation*}
Thus
\begin{equation}\label{2.26}
  \mathbb{Q}_{\delta_x}\langle f,I_t\rangle^2 = \frac{1}{\lambda^*}
  (\mathbb{V}ar_{\delta_x}\langle f,X_t\rangle-{\bf V}ar_{\delta_x}\langle f,\widetilde{X}_t\rangle).
\end{equation}

\begin{cor}\label{cor3.6}
Let $\{I_t\}_{t\geq0}$ be the process described in the Subsection 2.1 and $f\in\mathcal{P}$.
\begin{description}
\item{(i)} If $\alpha<2\gamma(f) b$, then
\begin{equation}\label{2.28}
\lim_{t\to\infty}  e^{-(\alpha/2)t}\mathbb{Q}_{\delta_x}\left(\langle f, I_t\rangle\right)=0,
\end{equation}
\begin{equation}\label{2.29}
\lim_{t\to\infty}e^{-\alpha t}\mathbb{V}_{\delta_x}\langle f,I_t\rangle=\frac{A}{\lambda^*}
 \int_0^{\infty} e^{\alpha s}\langle(T_sf)^2,\varphi\rangle\, ds=\frac{\sigma^2_f}{\lambda^*}.
\end{equation}
\item{(ii)} If $\alpha=2\gamma(f) b$, then
\begin{equation}\label{2.30}
 \lim_{t\to \infty}t^{-1/2} e^{-(\alpha/2)t}\mathbb{Q}_{\delta_x}\left(\langle f, I_t\rangle\right)=0,
 \end{equation}
and there exists $r\in\N$ such that
\begin{equation}\label{1.19}
  |t^{-1}e^{-\alpha t}\mathbb{V}_{\delta_x}\langle f,I_t\rangle|\lesssim 1+\|x\|^{2r}
\end{equation}
and
\begin{equation}\label{5.11}
\left|t^{-1}e^{-\alpha t}\mathbb{V}_{\delta_x}\langle f,I_t\rangle -\frac{A}{\lambda^*} \sum_{|p|=\gamma(f)}a_p^2\right|\lesssim t^{-1}(1+\|x\|^r),
\end{equation}
which in particular implies that
\begin{equation}\label{2.31}
\lim_{t\to\infty}t^{-1}e^{-\alpha t}\mathbb{V}_{\delta_x}\langle f,I_t\rangle=
  \frac{A}{\lambda^*} \sum_{|p|=\gamma(f)}a_p^2.
\end{equation}
\item{(iii)} If $\alpha>2\gamma(f) b$, then
\begin{equation}\label{3.17}
\lim_{t\to\infty}e^{-2(\alpha-\gamma(f) b)t}\mathbb{V}_{\delta_x}\langle f,I_t\rangle=
 \frac{\eta_f^2(x)}{\lambda^*}-\frac{1}{(\lambda^*)^2}(\sum_{|p|=\gamma(f)}a_p\phi_p(x))^2.
\end{equation}
\end{description}
\end{cor}
\textbf{Proof:}\quad Using \eqref{2.25} and Lemma~\ref{lem:2.2}, we can easily obtain the corollary. Here we just give the proof of \eqref{2.29}.
By \eqref{2.26}, we have
\begin{equation}\label{1.20}
  e^{-\alpha t}\mathbb{V}_{\delta_x}\langle f,I_t\rangle=\frac{1}{\lambda^*}e^{-\alpha t}\mathbb{V}ar_{\delta_x}\langle f,X_t\rangle -\frac{1}{\lambda^*}e^{-\alpha t}{\bf{V}}ar_{\delta_x}\langle f,\widetilde{X}_t\rangle-e^{-\alpha t}(\mathbb{Q}_{\delta_x}\langle f,I_t\rangle)^2.
\end{equation}
Using \eqref{24.14} and \eqref{2.15*}, we have
\begin{equation}\label{1.7}
{\bf V}ar_{\delta_x}\langle f, \widetilde{X}_t\rangle=(\psi^*)''(0+)e^{-\alpha^*t}\int_0^te^{-\alpha^*s}T_{t-s}[T_sf]^2(x)\,ds.
\end{equation}
By the fact that there exists $r\in\N$ such that $|T_tf(x)|\lesssim 1+\|x\|^r$,
we get $T_{t-s}[T_sf]^2(x)\lesssim (1+\|x\|^{2r})$. Thus
\begin{equation}\label{1.8}
 {\bf V}ar_{\delta_x}\langle f, \widetilde{X}_t\rangle\lesssim e^{-\alpha^*t}(1+\|x\|^{2r})\to 0, \quad t\to\infty.
\end{equation}
By \eqref{2.27}, $|\mathbb{Q}_{\delta_x}\langle f,I_t\rangle|\lesssim e^{\alpha t}|T_tf(x)|\lesssim  e^{(\alpha-\gamma(f)b)t}(1+\|x\|^{r})$, thus we have
\begin{equation}\label{1.21}
  \lim_{t\to\infty}e^{-\alpha t}(\mathbb{Q}_{\delta_x}\langle f,I_t\rangle)^2\lesssim \lim_{t\to\infty} e^{(\alpha-2\gamma(f)b)t}(1+\|x\|^{2r})=0.
\end{equation}
Now, using \eqref{2.2}, \eqref{1.8} and \eqref{1.21}, we easily get \eqref{2.29}. \hfill$\Box$

\begin{lemma}\label{4th-moment}
For $f\in\mathcal{P}$, it holds that
\begin{eqnarray}\label{4.15}
   \P_{\mu}(\langle f,\widetilde{X}_t\rangle-\P_{\mu}\langle f,\widetilde{X}_t\rangle)^4
   \lesssim \langle1+\|x\|^{4r},\mu\rangle+\langle1+\|x\|^{2r},\mu\rangle^2.
  \end{eqnarray}
  \end{lemma}
 \textbf{Proof:}\quad
By \eqref{2.40}, there exists $r\in \N$ such that $|T_tf(x)|\lesssim 1+\|x\|^r$. So by \eqref{4.17},
 $|(u_f^*)^{(1)}(x,t,0)|\lesssim 1+\|x\|^r$. By \eqref{1.8} and  \eqref{2.15*}, we have
  $|(u_f^*)^{(2)}(x,t,0)|\lesssim 1+\|x\|^{2r}$. Thus using \eqref{4.18}, we get
 $|(u_f^*)^{(3)}(x,t,0)|\lesssim 1+\|x\|^{3r}$.
Then by \eqref{4.19}, we have $|(u_f^*)^{(4)}(x,t,0)|\lesssim 1+\|x\|^{4r}$.
Now \eqref{4.15} follows immediately from \eqref{4.15*}. \hfill$\Box$

\section{Proofs of the main theorems}

In this section, we will prove the main results of this paper.
Recall that we assume that the initial measure $\mu$ is a finite measure on $\mathbb{R}^d$ with compact support, and that
$(X_t,\P_{\mu})$ and $(\Lambda_t,{\bf P}_{\mu})$ have the same law.
Thus in the remainder of this paper, we will replace $(X_t,\P_{\mu})$ by $(\Lambda_t,{\bf P}_{\mu})$.
Define
$${\cal L}_t=\{ u\in {\cal T}, \tau_u\le t<\sigma_u\},\quad t\ge 0.$$
From the construction of $\Lambda_t$, we have
\begin{equation}\label{3.6}
 \Lambda_{(t+s)}=\widetilde{X}^{t}_s+\sum_{u\in{\cal L}_t}I^{u,t}_s,
\end{equation}
where, conditioned on $\mathcal{G}_t$, $\widetilde{X}^{t}$ is a superprocess with the same law as $X$ under $\P^*_{\Lambda_t}$ and $I^{u,t}$ has the same law as $I$ under $\mathbb{Q}_{z_u(t)}$.
The processes $I^{u,t},u\in{\cal L}_t$, are independent.

\subsection{The large rate case: $\alpha>2b\gamma(f)$}

Recall that
$$
  H_t^{p}=e^{-(\alpha-|p|b)t}\langle\phi_p, X_t\rangle,\quad t\geq 0.
$$

\begin{lemma}\label{lem:1.2}
$H_t^p$ is a martingale under $\P_{\mu}$. Moreover, if $\alpha > 2|p|b$, we have $\sup_t\P_{\mu}(H_t^p)^2<\infty$, and therefore the limit
\begin{equation*}
  H_\infty^p:=\lim_{t\to \infty}H_t^p
\end{equation*}
exists $\P_{\mu}$-a.s. and in $L^2(\P_{\mu})$.
\end{lemma}

\textbf{Proof:}\quad
Since $\phi_p$ is an eigenfunction of $L$ corresponding to $-|p| b$,
by \eqref{2.17}, we have $\P_{\mu}H_t^p=\langle\phi_p,\mu\rangle$. Thus,
by the Markov property, we get that $H_t^p$ is a martingale.
Using \eqref{2.17} and \eqref{2.15}, we get
\begin{equation*}
\P_{\mu}\langle\phi_p, X_t\rangle^2=e^{2(\alpha-|p|b)t}\langle\phi_p,\mu\rangle^2+
Ae^{\alpha t}\int_{\mathbb{R}^d}\int_0^te^{(\alpha-2|p|b)s}T_{t-s}[\phi_{p}^2](x)\,ds\,\mu(dx).
\end{equation*}
Thus when $\alpha > 2|p|b$, we have by the definition of $H_t^p$,
\begin{eqnarray}
  \P_{\mu}(H_t^p)^2 &=& \langle\phi_{p},\mu\rangle^2+A\int_{\mathbb{R}^d}\int_0^t e^{-(\alpha-2|p|b)s}T_s[\phi_{p}^2](x)\,ds\,\mu(dx)\nonumber\\
   &\leq& \langle\phi_p,\mu\rangle^2+A\int_{\mathbb{R}^d}\int_0^\infty e^{-(\alpha-2|p|b)s}T_s[\phi_{p}^2](x)\,ds\,\mu(dx).\nonumber
\end{eqnarray}
Since $|\phi_{p}^2|\lesssim 1+\|x\|^{2|p|}$, by \eqref{e:ousmgp}, we have $|T_{s}[\phi_{p}^2](x)|\lesssim 1+\|x\|^{2|p|}$. Thus
\begin{equation}\label{7.14}
  \int_{\mathbb{R}^d}\int_0^\infty e^{-(\alpha-2|p|b)s}T_s[\phi_{p}^2](x)\,ds\,\mu(dx)\lesssim
  \int_{\mathbb{R}^d}(1+\|x\|^{2|p|})\,\mu(dx)<\infty,
\end{equation}
 from which the convergence asserted in the lemma follow easily.\hfill$\Box$

We now present the proof of  Theorem \ref{The:1.3}.

\textbf{Proof of Theorem \ref{The:1.3}:}\quad
 Define $M_t:=e^{-(\alpha-\gamma(f) b)t)}\langle \widetilde{f}, X_t\rangle$, where
 \begin{equation*}
   \widetilde{f}(x)=f(x)-\sum_{|p|=\gamma(f)}a_p\phi_p(x)=\sum_{n=\gamma(f)+1}^\infty\sum_{|p|=n}a_p\phi_p(x).
 \end{equation*}
It is clear that $\gamma(\widetilde{f})\geq\gamma(f)+1$. From Lemma~\ref{lem:2.2} and \eqref{2.40}, we have
\begin{description}
\item{(1)}
If $\alpha>2\gamma(\widetilde{f}))b$, then
\begin{eqnarray}
  \lim _{t\to\infty}e^{-2(\alpha-\gamma(\widetilde f)b)t)}\P_{\mu}\langle \widetilde{f}, X_t\rangle^2
\end{eqnarray}
exists, thus we have
\begin{eqnarray*}
  \P_{\mu}M_t^2 &=& e^{-2(\gamma(\widetilde f)-\gamma(f))b t}e^{-2(\alpha-\gamma(\widetilde f)b)t)} \P_{\mu}
\langle \widetilde{f}, X_t\rangle^2\\
   &=& O(e^{-2(\gamma(\widetilde f)-\gamma(f))b t})\to 0, \quad \mbox{ as } t\to\infty.
\end{eqnarray*}
\item{(2)}
If $\alpha=2\gamma(\widetilde{f})b$,
then $\lim _{t\to\infty}t^{-1}e^{-\alpha t}\P_{\mu}\langle \widetilde{f}, X_t\rangle^2$ exists.
Thus we have
\begin{eqnarray*}
   \P_{\mu}M_t^2 &=& te^{-2(\gamma(\widetilde{f})-\gamma(f)) t}(t^{-1}e^{-\alpha t} \P_{\mu}
   \langle \widetilde{f}, X_t\rangle^2)\\
   &=& O(te^{-2(\gamma(\widetilde{f})-\gamma(f)) t})\to 0, \quad \mbox{ as } t\to\infty.
\end{eqnarray*}
\item{(3)}
If $2\gamma(f) b<\alpha<2\gamma(\widetilde{f})b$,
then $\lim _{t\to\infty}e^{-\alpha t}\P_{\mu}\langle \widetilde{f}, X_t\rangle^2$ exists.
Thus we have
\begin{eqnarray*}
   \P_{\mu}M_t^2 &=& e^{-(\alpha-2\gamma(f) b) t}(e^{-\alpha t} \P_{\mu}
   \langle \widetilde{f}, X_t\rangle^2)\\
   &=& O(e^{-(\alpha-2\gamma(f) b) t})\to 0, \quad \mbox{ as } t\to\infty.
\end{eqnarray*}
\end{description}
Combining the three cases above, we get $\lim_{t\to\infty}M_t = 0$ in  $L^2(\P_{\mu})$.
Now using Lemma~\ref{lem:1.2}, we easily get the convergence in Theorem \ref{The:1.3}.\hfill$\Box$

\subsection{The small rate case: $\alpha<2\gamma(f)b$}

First, we recall some property of weak convergence.
For $f:\mathbb{R}^d\to\mathbb{R}$, let $\|f\|_L:=\sup_{x\ne y}|f(x)-f(y)|/\|x-y\|$ and
  $\|f\|_{BL}:=\|f\|_{\infty}+\|f\|_L$. For any distributions $\nu_1$ and $\nu_2$ on $\mathbb{R}^d$, define
\begin{equation*}
  \beta(\nu_1,\nu_2):=\sup\left\{\left|\int f\,d\nu_1-\int f\,d\nu_2\right|~:~\|f\|_{BL}\leq1\right\}.
\end{equation*}
Then $\beta$ is a metric. By \cite[Theorem 11.3.3]{Dudley}, the topology generated by this metric
is equivalent to the weak convergence topology.
 From the definition, we can easily see that, if $\nu_1$ and $\nu_2$ are the distributions of two $\R^d$-valued random variables $X$ and $Y$ respectively, then
\begin{equation}\label{5.20}
  \beta(\nu_1,\nu_2)\leq E\|X-Y\|\leq\sqrt{ E\|X-Y\|^2}.
\end{equation}

In the following, we will use the following elementary fact:
If X is a real-valued random variable with $E|X|^n<\infty$, then
\begin{equation}\label{3.20}
 \left| E(e^{i\theta X}-\sum_{m=0}^n\frac{(i\theta X)^m}{m!})\right|\leq
   \frac{|\theta|^n}{n!}E\left(|X|^n\left(\frac{|\theta X|}{n+1}\wedge 2\right)\right),
\end{equation}
which is an immediate consequence of the simple inequality
\begin{equation*}
  \left|e^{ix}-\sum_{m=0}^n\frac{(ix)^m}{m!}\right|\leq \min\left(\frac{|x|^{n+1}}{(n+1)!}, \frac{2|x|^n}{n!}\right).
\end{equation*}

\smallskip

Now we are ready to prove Theorem \ref{The:1.4}.

\textbf{Proof of Theorem \ref{The:1.4}:}\quad
 We define an ${\mathbb R}^2$-valued random variable $U_1(t)$ by
\begin{equation}\label{5.16}
   U_1(t):=(e^{-\alpha t}\left\|\Lambda_t\right\|, e^{-(\alpha/2) t}\langle f,\Lambda_t\rangle).
\end{equation}
To get the conclusion of Theorem \ref{The:1.4}, it suffices to show that under ${\bf P}_\mu$,
\begin{equation}\label{5.17}
   U_1(t)\stackrel{d}{\to}(W_\infty, \sqrt{W_\infty}G_1(f)),
\end{equation}
where $G_1(f)\sim\mathcal{N}(0,\sigma_f^2)$.
Let $s, t >0$ and write
\begin{equation*}
  U_1(s+t)=(e^{-\alpha (s+t)}\left\|\Lambda_{s+t}\right\|, e^{-(\alpha/2) (s+t)}\langle f,\Lambda_{s+t}\rangle).
\end{equation*}
Recall the representation \eqref{3.6}.
Define
\begin{equation}\label{e:new}
Y_s^{u,t}:=e^{-\alpha s/2}\langle f,I_{s}^{u,t}\rangle\quad\mbox{and}\quad y_s^{u,t}:={\bf P}_{\mu}(Y^{u,t}_s|\mathcal{G}_t).
\end{equation}
$Y_s^{u,t}$ has the same law as $Y_s:=e^{-\alpha s/2}\langle f,I_{s}\rangle$ under $\mathbb{Q}_{\delta_{Z_u(t)}}$.
Then we have
\begin{eqnarray}\label{1.24}
     &&e^{-(\alpha/2) (s+t)}\langle f,\Lambda_{s+t}\rangle\nonumber\\
  &=& e^{-(\alpha/2) (s+t)} \langle f,\widetilde{X}_{s}^t\rangle+e^{-(\alpha/2) t}\sum_{u\in\cL_t}Y_s^{u,t}\nonumber\\
  &=& e^{-(\alpha/2) (s+t)}(\langle f,\widetilde{X}_{s}^t\rangle-{\bf P}_{\mu}(\langle f,\widetilde{X}_{s}^t\rangle|\mathcal{G}_t)\nonumber\\
  &&+e^{-(\alpha/2) t}\sum_{u\in\cL_t}(Y_s^{u,t}-y^{u,t}_s)+ e^{-(\alpha/2) (t+s)}{\bf P}_{\mu}(\langle f,\Lambda_{s+t}\rangle|\mathcal{G}_t)\nonumber\\
  &=:& J_0(s,t)+J_1(s,t)+J_2(s,t).
\end{eqnarray}
Put $\widetilde{V}_s(x):={\bf V}ar_{\delta_x}\langle f,\widetilde{X}_{s}\rangle$. Then
$$
{\bf P}_{\mu}J_0(s,t)^2=e^{-\alpha(t+s)}{\bf P}_{\mu}\langle \widetilde{V}_s,\Lambda_t\rangle=e^{-\alpha s}\langle T_t\widetilde{V}_s,\mu\rangle.
$$
By \eqref{1.8}, there exists $r\in \N$ such that $\widetilde{V}_s(x)\lesssim e^{-\alpha^*s}(1+\|x\|^{2r})$.
Thus \begin{equation}\label{J1to0}{\bf P}_{\mu}J_0(s,t)^2\lesssim e^{-\alpha s}e^{-\alpha^*s}\int_{\mathbb{R}^d}(1+\|x\|^{2r})\mu(dx).\end{equation}

Next we consider $J_{2}(s,t)$.
By the Markov property and \eqref{2.17}, we have
\begin{equation*}
  J_{2}(s,t)=e^{-(\alpha/2) (s+t)}e^{\alpha s}\langle T_{s}f, \Lambda_t\rangle.
\end{equation*}
Thus, by \eqref{2.14} and \eqref{2.17}, we have
\begin{eqnarray}
  {\bf P}_{\mu}J_{2}(s,t)^2 &=& A e^{\alpha s}\int_{\mathbb{R}^d}\int_0^t
   e^{\alpha u}T_{t-u}[T_{s+u}f]^2(x)\,du\,\mu(dx) +e^{\alpha(t+s)}\langle T_{t+s}f, \mu\rangle^2\nonumber\\
  &\lesssim& e^{(\alpha-2\gamma(f)b)s}\int_{\mathbb{R}^d}(1+\|x\|^{2r})\,\mu(dx),\label{10.1}
\end{eqnarray}
here the last inequality follows from the fact that there exists $r\in\N$
such that
$$
|T_{s+u}f|(x)\lesssim e^{-\gamma(f)b(u+s)}(1+\|x\|^r).
$$
Thus by \eqref{J1to0} and \eqref{10.1}, we have
\begin{equation}\label{10.2}
  \lim_{s\to\infty}\limsup_{t\to\infty}{\bf P}_{\mu}(J_0(s,t)+J_2(s,t))^2=0.
\end{equation}

Now we consider $J_1(s,t)$.
We define an ${\mathbb R}^2$-valued random variable $U_2(s,t)$ by
\begin{eqnarray*}
U_2(s,t):=\left(e^{-\alpha t}\|\Lambda_t\|,J_1(s,t) \right).
\end{eqnarray*}
We claim that under ${\bf P}_\mu$,
\begin{equation}\label{10.5}
  U_2(s,t)\stackrel{d}{\to}(W_\infty, \sqrt{W_\infty}G_1(s)), \quad \mbox{ as } t\to\infty,
\end{equation}
where $G_1(s)\sim\mathcal{N}(0,\sigma^2_f(s))$ and $\sigma^2_f(s)$ will be given later.
Denote the characteristic function of $U_2(s,t)$ under ${\bf P}_\mu$ by
$\kappa(\theta_1,\theta_2,s,t)$:
\begin{eqnarray}
  \kappa(\theta_1,\theta_2,s,t)
  &=&{\bf P}_{\mu}\left(\exp\left\{i\theta_1e^{-\alpha t}\|\Lambda_t\|+i\theta_2e^{-(\alpha/2) t}\sum_{u\in\cL_t}(Y_s^{u,t}-y_s^{u,t})\right\}\right)\nonumber\\
  &=& {\bf P}_{\mu}\left(\exp\{i\theta_1e^{-\alpha t}\|\Lambda_t\|\}\prod_{u\in\cL_t}h_s(z_u(t),e^{-(\alpha/2) t}\theta_2)\right)\nonumber\\
  &=& {\bf P}_{\mu}\left(\exp\{i\theta_1e^{-\alpha t}\|\Lambda_t\|\}\exp\left\{\lambda^*\langle h_s(\cdot,~e^{-(\alpha/2)t}\theta_2)-1,\Lambda_t\rangle\right\}\right),\label{10.9}
\end{eqnarray}
where $h_s(x,\theta)=\mathbb{Q}_{\delta_x}e^{i\theta(Y_s-\mathbb{Q}_{\delta_x}Y_s)}$.
The last equality in the display above follows from the fact that given $\Lambda_t$,
$Z_t$ is a Poisson random measure  with intensity $\lambda^*\Lambda_t$.
Define
\begin{equation*}
   e_s(x,\theta):=h_s(x,\theta)-1+\frac{1}{2}\theta^2\mathbb{V}_{\delta_x}Y_s
\end{equation*}
and $V_s(x):=\mathbb{V}_{\delta_x}Y_s$. Then
\begin{eqnarray*}
   \exp\left\{\lambda^*\langle h_s(\cdot,~e^{-(\alpha/2)t}\theta_2)-1,\Lambda_t\rangle\right\}&=& \exp\left\{-\lambda^*\frac{1}{2}\theta_2^2e^{-\alpha t}\langle V_s, \Lambda_t\rangle\right\} \exp\left\{\lambda^*\langle e_s(\cdot,e^{-(\alpha/2)t}\theta_2), \Lambda_t\rangle\right\}\\
   &=&J_{1,1}(s,t)J_{1,2}(s,t).
\end{eqnarray*}
By \eqref{3.20}, we have
$$
|e_s(x,e^{-(\alpha/2)t}\theta_2)|\leq \theta_2^2e^{-\alpha t}\mathbb{Q}_{\delta_x}\left(|Y_s-\mathbb{Q}_{\delta_x}Y_s|^2(\frac{e^{-(\alpha/2)t}\theta_2|Y_s-\mathbb{Q}_{\delta_x}Y_s|}{6}\wedge 1)\right).
$$
Let
$$
g(x,s,t):=\mathbb{Q}_{\delta_x}\left(|Y_s-\mathbb{Q}_{\delta_x}Y_s|^2(
\frac{e^{-(\alpha/2)t}\theta_2|Y_s-\mathbb{Q}_{\delta_x}Y_s|}{6}\wedge 1)\right).
$$
We notice $g(x,s,t)\downarrow0$ as  $t\uparrow\infty$. By \eqref{2.17},
\begin{equation*}
  {\bf P}_{\mu}|\langle e_s(\cdot,e^{-(\alpha/2)t}\theta_3), \Lambda_t\rangle|\leq \theta^2_2 ~\langle~ T_t(g(\cdot,s,t)),~\mu~\rangle\to 0, \quad \mbox{as}\quad t\to\infty.
\end{equation*}
Thus $\lim_{t\to\infty}\langle e_s(\cdot,e^{-(\alpha/2)t}\theta_2), \Lambda_t\rangle=0$ in probability,
which implies $\lim_{t\to\infty}J_{1,2}(s,t)=1$  in probability.
Furthermore, by Remark \ref{rem:large}, we have
\begin{equation*}
\lim_{t\to\infty}e^{-\alpha t}\langle V_s, \Lambda_t\rangle
=\langle V_s,~\varphi\rangle W_\infty\quad \mbox{in probability},
\end{equation*}
which implies that $\lim_{t\to\infty}J_{1,1}(s,t)=\exp\left\{-\frac{1}{2}\theta_2^2\sigma_f^2(s) W_\infty\right\}$,
where $\sigma_f^2(s):= \lambda^*\langle V_s,~\varphi\rangle$. Thus
\begin{equation}\label{10.8}
 \lim_{t\to\infty}\exp\left\{\lambda^*\langle h_s(\cdot,~e^{-(\alpha/2)t}\theta_2)-1,\Lambda_t\rangle\right\}=
 \exp\left\{-\frac{1}{2}\theta_2^2\sigma_f^2(s) W_\infty\right\} \quad \mbox{in probability}.
\end{equation}
Since $h_s(x,\theta)$ is a characteristic function, its real part is less than 1, which implies
$$
\left|\exp\left\{\lambda^*\langle h_s(\cdot,~e^{-(\alpha/2)t}\theta_2)-1,\Lambda_t\rangle\right\}\right|\leq 1.
$$
Hence by the dominated convergence theorem, we get
\begin{equation}\label{10.4}
 \lim_{t\to\infty}\kappa(\theta_1,\theta_2,s,t)= {\bf P}_{\mu}\exp\left\{i\theta_1W_\infty\right\}
  \exp\left\{-\frac{1}{2}\theta_2^2\sigma_f^2(s) W_\infty\right\},
\end{equation}
which implies our claim \eqref{10.5}.
Thus, we easily get that under ${\bf P}_\mu$,
\begin{eqnarray*}
U_3(s,t):=\left(e^{-\alpha (t+s)}\|\Lambda_{t+s}\|,J_1(s,t) \right)\stackrel{d}{\to}(W_\infty, \sqrt{W_\infty}G_1(s)), \quad \mbox{ as } t\to\infty.
\end{eqnarray*}
By \eqref{2.29}, we have $\lim_{s\to\infty}V_s(x)= \frac{\sigma^2_f}{\lambda^*}$,
thus $\lim_{s\to\infty}\sigma_f^2(s)=\sigma^2_f$. So
\begin{equation}\label{10.7}
\lim_{s\to\infty}\beta(G_1(s),G_1(f))=0.
\end{equation}
Let $\mathcal{L}(s+t)$ and $\widetilde{\mathcal{L}}(s,t)$ be the distributions of $U_1(s+t)$ and $U_3(s,t)$
respectively, and let $\mathcal{L}(s)$ and $\mathcal{L}$ be the distributions of $(W_\infty, \sqrt{W_\infty}G_1(s))$
and $(W_\infty, \sqrt{W_\infty}G_1(f))$ respectively. Then, using \eqref{5.20}, we have
\begin{eqnarray}\label{10.11}
  \limsup_{t\to\infty}\beta(\mathcal{L}(s+t),\mathcal{L})&\leq&
  \limsup_{t\to\infty}[\beta(\mathcal{L}(s+t),\widetilde{\mathcal{L}}(s,t))+\beta(\widetilde{\mathcal{L}}(s,t),\mathcal{L}(s))+\beta(\mathcal{L}(s),\mathcal{L})]\nonumber\\
 &\leq &\limsup_{t\to\infty}({\bf P}_{\mu}(J_0(s,t)+J_2(s,t))^2)^{1/2}+0+\beta(\mathcal{L}(s),\mathcal{L}).
\end{eqnarray}
Using this and the definition of $\limsup_{t\to\infty}$, we easily get that
$$
\limsup_{t\to\infty}\beta (\mathcal{L}(t),\mathcal{L})=\limsup_{t\to\infty}\beta (\mathcal{L}(s+t),\mathcal{L})
\le \limsup_{t\to\infty}({\bf P}_{\mu}(J_0(s,t)+J_2(s,t))^2)^{1/2}+\beta(\mathcal{L}(s),\mathcal{L}).
$$
Letting $s\to\infty$, we get $ \limsup_{t\to\infty}\beta(\mathcal{L}(t),\mathcal{L})=0$.
The proof is now complete. \hfill$\Box$

\subsection{Proof of Theorem \ref{The:2.1}}
In this section we consider the case $\alpha>2\gamma(f)b$ and $f_{(c)}= 0$.
Recall the decomposition of $\Lambda_t$ under ${\bf P}_{\delta_x}$ in \eqref{2.25},
we have for $|p|=m<\alpha/(2b)$,
\begin{equation}\label{8.27}
  H_s^p=e^{-(\alpha-mb)s}\langle\phi_p, \widetilde{X}_s\rangle+\sum_{j=1}^Ne^{-(\alpha-mb)s}\langle\phi_p, I_s^j\rangle.
\end{equation}
Let $$\widetilde{H}^p_s:=e^{-(\alpha-mb)(s)}\langle\phi_p,I_s\rangle. $$
Then under ${\bf P}_{\delta_x}$, the processes $\{e^{-(\alpha-mb)s}\langle\phi_p,I_s^{j}\rangle ,s\geq 0\},
~j=1,2\dots$  are i.i.d. with a common law equal to that of
$\{\widetilde{H}^p_s,s\geq 0\}$ under $\mathbb{Q}_{\delta_x}$.
Since $\phi_p$ is an eigenvalue of $L$ corresponding to $-|p|b$, we have
\begin{equation}\label{8.41}
  {\bf P}_{\delta_x}\langle\phi_p, \widetilde{X}_s\rangle=e^{-(\alpha^*+mb)s}\phi_p(x)\to 0,\quad \mbox{ as } s\to\infty.
\end{equation}
Thus, by \eqref{1.8}, we have that as $s\to\infty$,
\begin{equation}\label{8.45}
  {\bf P}_{\delta_x}(\langle\phi_p, \widetilde{X}_s\rangle)^2\lesssim e^{-\alpha^*s}(1+\|x\|^{2|p|})\to 0,
\end{equation}
which implies $\lim_{s\to\infty}e^{-(\alpha-mb)s}\langle\phi_p, \widetilde{X}_s\rangle=0$ in $L^2({\bf P}_{\delta_x})$.
By Lemma \ref{lem:1.2}, $\lim_{s\to\infty}H_s^p=H^p_\infty$ in $L^2({\bf P}_{\delta_x})$. Thus
\begin{equation}\label{8.28}
\lim_{s\to\infty}\sum_{j=1}^Ne^{-(\alpha-mb)s}\langle\phi_p, I_s^j\rangle = H^p_\infty
\quad \mbox{ in } L^2({\bf P}_{\delta_x}).
\end{equation}
From the fact that $N$ is independent of $I^{j}$, we have for any $s,t\geq 0$,
\begin{eqnarray*}
  &&{\bf P}_{\delta_x}\left[\sum_{j=1}^N(e^{-(\alpha-mb)s}\langle\phi_p, I_s^j\rangle-e^{-(\alpha-mb)t}\langle\phi_p, I_t^j\rangle)\right]^2\\
  &\geq& {\bf P}_{\delta_x}[(e^{-(\alpha-mb)s}\langle\phi_p, I_s\rangle-e^{-(\alpha-mb)t}\langle\phi_p, I_t\rangle)^2;N=1]\\
  &=&{\bf P}_{\delta_x}(N=1)\mathbb{Q}_{\delta_x}(\widetilde{H}^p_s-\widetilde{H}^p_t)^2.
\end{eqnarray*}
By \eqref{8.28}, we get for any $x\in{\mathbb R}^d$,
\begin{equation}\label{8.29}
  \mathbb{Q}_{\delta_x}(\widetilde{H}^p_s-\widetilde{H}^p_t)^2\to 0,\quad s,t\to\infty.
\end{equation}
 Thus $\widetilde{H}^p_s$ converges in $L^2(\mathbb{Q}_{\delta_x})$. Let
$$\widetilde{H}^p_\infty:=\lim_{s\to\infty}\widetilde{H}^p_s \quad\mbox{ in }L^2(\mathbb{Q}_{\delta_x}),$$
which implies $H^{j, p}_\infty:=\lim_{s\to\infty}\langle\phi_p,I_s^{j}\rangle e^{-(\alpha-mb)s}$
exists in $L^2({\bf P}_{\delta_x})$.
Furthermore, $H^{j, p}_\infty$, under ${\bf P}_{\delta_x}$, are i.i.d.
with a common law equal to that of  $\widetilde{H}^p_\infty$ under $\mathbb{Q}_{\delta_x}$.
Hence by \eqref{8.28}, it is easy to get
\begin{equation}\label{8.30}
  H_\infty^p=\sum_{j=1}^NH^{j, p}_\infty,\quad {\bf P}_{\delta_x}\mbox{- a.s.}
\end{equation}

Recall the decomposition of $\Lambda_{t+s}$ in \eqref{3.6}.
By Lemma \ref{lem:1.2}, we have for $|p|=m<\alpha/(2b)$,
\begin{eqnarray}\label{8.2}
  H^p_{t+s}=e^{-(\alpha-mb)(s+t)}\langle\phi_p,\widetilde{X}^t_s\rangle+e^{-(\alpha-mb)t}
  \sum_{u\in{\cL}_t}e^{-(\alpha-mb)s}\langle\phi_p,I_s^{u,t}\rangle .
\end{eqnarray}
From the definition of $\widetilde{X}^t_s$, using \eqref{1.8} and \eqref{8.41}, we have
\begin{eqnarray*}
  {\bf P}_{\mu}(\langle\phi_p,\widetilde{X}^t_s\rangle)^2 &\leq & 2{\bf P}_{\mu}(\langle\phi_p,\widetilde{X}^t_s\rangle-{\bf P}_{\mu}(\langle\phi_p,\widetilde{X}^t_s\rangle|\mathcal{F}_t))^2 +2{\bf P}_{\mu}\left({\bf P}_{\mu}(\langle\phi_p,\widetilde{X}^t_s\rangle|\mathcal{F}_t)\right)^2\\
  &=& 2{\bf P}_{\mu}\langle {\bf V}ar_{\delta_\cdot}\langle \phi_p, \widetilde{X}_s\rangle, \Lambda_t\rangle+2{\bf P}_{\mu}\langle{\bf P}_{\delta_\cdot}\langle\phi_p,\widetilde{X}_s\rangle,\Lambda_t \rangle^2\to 0,\quad \mbox{ as }s\to\infty.
\end{eqnarray*}
Hence $\lim_{s\to\infty}e^{-(\alpha-mb)(s+t)}\langle\phi_p,\widetilde{X}^t_s\rangle=0$ in $L^2({\bf P}_{\mu})$. Thus
$\lim_{s\to\infty}e^{-(\alpha-mb)(s+t)}\langle\phi_p,\widetilde{X}^t_s\rangle=0$ in $L^2({\bf P}_{\mu})$. Thus
\begin{equation}\label{8.42}
\lim_{s\to\infty}e^{-(\alpha-mb)t}\sum_{u\in{\cL}_t}\langle\phi_p,I_s^{u,t}\rangle e^{-(\alpha-mb)s}=
  H_\infty^p \quad\mbox{ in }L^2({\bf P}_{\mu}).
\end{equation}
Note that under ${\bf P}_{\mu}$, given $Z_t$, $e^{-(\alpha-mb)(s)}\langle\phi_p,I_s^{u,t}\rangle $
has the same law as $\widetilde{H}^p_s$
under $\mathbb{Q}_{\delta_{Z_u(t)}}$. Thus by \eqref{8.29},
for each $u\in {\cL}_t$, $e^{-(\alpha-mb)s}\langle\phi_p,I_s^{u,t}\rangle $ converges in $L^2({\bf P}_{\mu})$
to a limit, denoted as $H^{u,t, p}_\infty$. Furthermore, given $Z_t$, $H^{u,t,p}_\infty$ has the same law
as $\widetilde{H}^p_\infty$ under $\mathbb{Q}_{\delta_{Z_u(t)}}$.

We claim that, for each $t\geq 0$,
\begin{equation}\label{8.4}
   H_\infty^p=e^{-(\alpha-mb)t}\sum_{u\in{\cL}_t}H^{u,t, p}_\infty.
\end{equation}
In fact,
\begin{eqnarray*}
  {\bf P}_{\mu}(\sum_{u\in{\cL}_t}e^{-(\alpha-mb)s}\langle\phi_p,I_s^{u,t}\rangle -H^{u,t, p}_\infty)^2
  &\leq&  {\bf P}_{\mu}|Z_t|\sum_{u\in{\cL}_t}(e^{-(\alpha-mb)s}\langle\phi_p,I_s^{u,t}\rangle -H^{u,t, p}_\infty)^2  \\
   &=&  {\bf P}_{\mu}|Z_t|\sum_{u\in{\cL}_t}\mathbb{Q}_{\delta_{Z_u(t)}}(\widetilde{H}^{p}_s-\widetilde{H}^p_\infty)^2.
\end{eqnarray*}
By \eqref{2.26}, we have
\begin{eqnarray*}
   \mathbb{Q}_{\delta_x}(\widetilde{H}^p_s)^2 \leq  \frac{1}{\lambda^*}\mathbb{V}ar_{\delta_x}({H}^p_s)\leq \frac{1}{\lambda^*}\mathbb{P}_{\delta_x}(H_s^p)^2\lesssim 1+\|x\|^{2|p|}.
\end{eqnarray*}
Thus $\mathbb{Q}_{\delta_{x}}(\widetilde{H}^{p}_s-\widetilde{H}^p_\infty)^2\leq 2\sup_{s\geq 0}
 \mathbb{Q}_{\delta_x}(\widetilde{H}^p_s)^2 \lesssim 1+\|x\|^{2|p|}$. We can easily get that
$$
{\bf P}_{\mu}|Z_t|\langle(1+\|\cdot\|^{2|p|}),Z_t\rangle< \infty.
$$
So by the dominated convergence theorem, we have $\lim_{s\to\infty} {\bf P}_{\mu}(\sum_{u\in{\cL}_t}
e^{-(\alpha-mb)s}\langle\phi_p,I_s^{u,t}\rangle-H^{u,t, p}_\infty)^2= 0$.
Now the claim \eqref{8.4} follows easily from \eqref{8.42}.

Define
 $$
 H^{u,t}_\infty:=\sum_{\gamma(f)\leq m<\alpha/2b}\sum_{|p|=m}a_pH^{u,t, p}_\infty
 \quad \mbox{and}\quad \widetilde{H}_\infty:=\sum_{\gamma(f)\leq m<\alpha/2b}\sum_{|p|=m}a_p\widetilde{H}^p_\infty.
 $$
 Recall the definition of $H_\infty$ in \eqref{8.47}. By \eqref{8.30},
 we have $$H_\infty=\sum_{u\in{\cL}_0}H^{u,0}_\infty.$$
 Under ${\bf P}_{\delta_x}$, $H^{u,0}_\infty$ are i.i.d. with a common law equal to that of
 $\widetilde{H}_\infty$ under $\mathbb{Q}_{\delta_x}$. Thus we have
\begin{eqnarray}
  {\bf P}_{\delta_x}H_\infty &=& \lambda^*\mathbb{Q}_{\delta_x}\widetilde{H}_\infty, \label{8.19}\\
  {\bf V}ar_{\delta_x}H_\infty &=& \lambda^*\mathbb{Q}_{\delta_x}(\widetilde{H}_\infty)^2.\label{8.20}
\end{eqnarray}
On the other hand, by Lemma \ref{lem:1.2}, we get
$$
\lim_{t\to\infty}\sum_{\gamma(f)\le m<\alpha/2b}e^{-(\alpha-mb)t}\sum_{|p|=m}a_p\langle\phi_p,\Lambda_t\rangle= H_\infty,
\quad \mbox{in} ~L^2({\bf P}_{\delta_x}).
$$
It follows that
\begin{equation}\label{8.21}
  {\bf P}_{\delta_x}H_\infty=f_{(s)}(x),
\end{equation}
and by \eqref{2.24},
\begin{equation}\label{8.22}
  {\bf V}ar_{\delta_x}H_\infty=A\int_0^\infty e^{-\alpha s}T_s
    \left(\sum_{\gamma(f)\le m<\alpha/2b}e^{mbs}\sum_{|p|=m}a_p\phi_p\right)^2(x)\,ds.
\end{equation}

\smallskip

 \textbf{Proof of Theorem \ref{The:2.1}:}\quad
 By \eqref{8.4}, we have
 $$
 \sum_{\gamma(f)\le m<\alpha/2b}e^{(\alpha-mb)t}\sum_{|p|=m}a_pH^p_\infty=\sum_{u\in\cL_t}H^{u,t}_\infty.
 $$
 Consider the $\mathbb{R}^2$-valued random variable $U_1(t)$:
 \begin{equation}\label{8.5}
   U_1(t):= \left(e^{-\alpha t}\|\Lambda_t\|,~e^{-(\alpha/2)t}(\langle
   f,\Lambda_t\rangle-\sum_{u\in{\cL}_t}H^{u,t}_\infty)\right).
 \end{equation}
 To get the conclusion of Theorem \ref{The:2.1}, it suffices to show that
\begin{equation}\label{8.18}
 U_1(t)\stackrel{d}{\to}(W_\infty, \sqrt{W_\infty}G_3(f)).
\end{equation}
 Denote the characteristic function of $U_1(t)$  with respect to  ${\bf P}_{\mu}$ by $\kappa_1(\theta_1,\theta_2,t)$ and let $h(x,\theta):=\mathbb{Q}_{\delta_x}\exp\{i\theta\widetilde{H}_\infty\}$. Then we have
 \begin{eqnarray}\label{8.6}
  &&\kappa_1(\theta_1,\theta_2,t)\nonumber\\
  &=&
   {\bf P}_{\mu}\exp\left\{i\theta_1e^{-\alpha t}\|\Lambda_t\|+i\theta_2e^{-(\alpha/2)t}(\langle
   f,\Lambda_t\rangle-\sum_{u\in{\cL}_t}H^{u,t}_\infty)\right\}\nonumber\\
   &=&{\bf P}_{\mu}\exp\left\{i\theta_1e^{-\alpha t}\|\Lambda_t\|\right\}\exp\left\{i\theta_2e^{-(\alpha/2)t}\langle f,\Lambda_t\rangle\right\}\prod_{u\in{\cL}_t}h\left(Z_u(t),-\theta_2e^{-(\alpha/2)t}\right)\nonumber\\
   &=&{\bf P}_{\mu}\exp\left\{i\theta_1e^{-\alpha t}\|\Lambda_t\|\right\}\exp\left\{i\theta_2e^{-(\alpha/2)t}\langle f,\Lambda_t\rangle+\lambda^*\langle h(\cdot,-\theta_2e^{-(\alpha/2)t})-1, \Lambda_t\rangle\right\}.
 \end{eqnarray}
The third equality above follows from the fact that, given $\Lambda_{t}$,
$Z_{t}$ is a Poisson point process with density $\lambda^*\Lambda_{t}$. By \eqref{8.19} and \eqref{8.21}, we get $\mathbb{Q}_{\delta_x}\widetilde{H}_\infty=f_{(s)}(x)/\lambda^*$.
Let
$$
e(x,\theta):=h(x,\theta)-1-\frac{i\theta}{\lambda^*}f_{(s)}(x) +\frac{1}{2}
\mathbb{Q}_{\delta_x}(\widetilde{H}_\infty )^2\theta^2
$$
and $V(x):={\bf V}ar_{\delta_x}H_\infty$. Then, by \eqref{8.20}, we have
 \begin{eqnarray*}
    &&i\theta_2e^{-(\alpha/2)t}\langle f,\Lambda_t\rangle+\lambda^*\langle h(\cdot,-\theta_2e^{-(\alpha/2)t})-1, \Lambda_t\rangle \\
   &=& i\theta_2e^{-(\alpha/2)t}\langle f_{(l)},\Lambda_t\rangle-\frac{1}{2}\theta_2^2e^{-\alpha t}\langle V,\Lambda_t\rangle+\lambda^*\langle e(\cdot,-\theta_2e^{-(\alpha/2)t}),\Lambda_t\rangle \\
    &=:& J_1(t)+J_2(t)+J_3(t).
 \end{eqnarray*}
By \eqref{3.20}, we have
\begin{equation}\label{8.9}
 |e(x,\theta)|\leq \theta^2\mathbb{Q}_{\delta_x}\left(|\widetilde{H}_\infty|^2\left(\frac{\theta|\widetilde{H}_\infty|}{6}\wedge 1\right)\right),
\end{equation}
which implies that
\begin{eqnarray*}
  |J_3(t)|&\leq& \theta_2^2e^{-\alpha t}\langle g(\cdot,~t),\Lambda_t\rangle,
\end{eqnarray*}
where
$$
g(x,t):=\mathbb{Q}_{\delta_x}\left(|\widetilde{H}_\infty|^2\left(\frac{e^{-(\alpha/2)t}\theta_2|\widetilde{H}_\infty|}{6}
\wedge 1\right)\right).
$$
It is clear that $g(x,t)\downarrow 0$ as  $t\uparrow\infty$.
Thus
\begin{equation}
  {\bf P}_{\mu}|J_3(t)|\leq \theta^2_2 ~\langle~ T_t(g(\cdot,t)),~\mu~\rangle\to 0, \quad \mbox{as}\quad t\to\infty,
\end{equation}
which implies $\lim_{t\to\infty}J_3(t)= 0$ in probability.
By Remark \ref{rem:large}, we have
\begin{equation}\label{8.13}
\lim_{t\to\infty} e^{-\alpha t}\langle V , \Lambda_t\rangle= \langle V,~\varphi\rangle W_\infty \quad \mbox{in probability}.
\end{equation}
Recall that $\lim_{t\to\infty}e^{-\alpha t}\|\Lambda_t\|=W_\infty,$ ${\bf P}_\mu$-a.s. Therefore
\begin{equation}\label{8.23}
  \lim_{t\to\infty}\exp\left\{i\theta_1e^{-\alpha t}\|\Lambda_t\|\right\}\exp\{J_2(t)+J_3(t)\}=
  \exp\{i\theta_1W_\infty\}\exp\{-\frac{1}{2}\theta_2^2\langle V,~\varphi\rangle W_\infty\}\quad \mbox{ in probability}.
\end{equation}
Thus by the dominated convergence theorem, we get that as $t\to\infty$,
\begin{equation}\label{8.14}
  \left|\kappa_1(\theta_1,\theta_2,t)
  -{\bf P}_{\mu}\exp\left\{i\theta_2e^{-(\alpha/2)t}\langle f_{(l)},\Lambda_t\rangle\right\}\exp\{i\theta_1W_\infty\}\exp\{-\frac{1}{2}\theta_2^2\langle V,~\varphi\rangle W_\infty\}\right|\to 0.
\end{equation}
Since $\alpha<2\gamma(f_{(l)})b$, by Theorem \ref{The:1.4}, we have that as $t\to\infty$,
\begin{equation}\label{8.25}
  (e^{-\alpha t}\|\Lambda_t\|,e^{-(\alpha/2)t}\langle f_{(l)},\Lambda_t\rangle)\stackrel{d}{\to}(W_\infty,\sqrt{W_\infty}G_1(f_{(l)})),
\end{equation}
where $G_1(f_{(l)})\sim\mathcal{N}(0,\sigma^2_{f_{(l)}})$.
 Therefore, 
 \begin{eqnarray}\label{8.26}
  &&\lim_{t\to\infty} {\bf P}_{\mu}\exp\left\{i\theta_2e^{-(\alpha/2)t}\langle f_{(l)},\Lambda_t\rangle\right\}e^{i\theta_1W_\infty}\exp\{-\frac{1}{2}\theta_2^2\langle V,~\varphi\rangle W_\infty\}\nonumber\\ 
  &=&{\bf P}_{\mu}e^{i\theta_1W_\infty}\exp\{-\frac{1}{2}\theta_2^2(\sigma^2_{f_{(l)}}+\langle V,~\varphi\rangle )W_\infty\}.
\end{eqnarray}
By \eqref{8.22}, we get
\begin{equation*}
   \langle V,~\varphi\rangle =A\sum_{\gamma(f)\le m<\alpha/2b}\frac{1}{\alpha-2mb}\sum_{|p|=m}a_p^2.
\end{equation*}
The proof is now complete.\hfill$\Box$

\subsection{The  critical case: $\alpha=2\gamma(f)b$}

To prove Theorem \ref{The:1.5}, we need the following lemma. The idea of the proof is mainly from \cite{Ath71}.
\begin{lemma}\label{lem:5.5}
Assume $f\in \mathcal{P}$ satisfies $\alpha = 2\gamma(f)b$.
Define $T^\alpha_tf(x):=e^{\alpha t}T_tf(x)=\mathbb{P}_{\delta_x}\langle f,X_t\rangle$ and
$$
S_tf:=t^{-1/2}e^{-(\alpha/2) t}(\langle f,X_t\rangle-{T}^\alpha_tf(x)).
$$
Then for any $c>0$ and $\delta>0$, we have
\begin{equation}\label{4.1}
\lim_{t\to\infty}\mathbb{P}_{\delta_x}\left(|S_tf|^2;|S_tf|>ce^{\delta t }\right)=0.
\end{equation}
\end{lemma}

\textbf{Proof:}\quad
We write $t=[t]+\epsilon_t$, where $[t]$ is the integer part of $t$. Then
\begin{eqnarray}\label{4.35}
  S_{t}f&=&t^{-1/2}e^{-(\alpha/2)t}\left(\langle f,X_t\rangle-\langle T^\alpha_{\epsilon_t}f, X_{[t]}\rangle\right)+t^{-1/2}e^{-(\alpha/2)t}\left(\langle T^\alpha_{\epsilon_t}f, X_{[t]}\rangle-T^\alpha_{t}f(x)\right)\nonumber\\
  &=&t^{-1/2}e^{-(\alpha/2)t}\left(\langle f,X_t\rangle-\langle T^\alpha_{\epsilon_t}f,X_{[t]}\rangle\right)
  +\left(\frac{[t]}{t}\right)^{1/2}e^{-\alpha\epsilon_t/2}S_{[t]}(T^\alpha_{\epsilon_t}f).
\end{eqnarray}
Thus
\begin{eqnarray*}
  &&\mathbb{P}_{\delta_x}\left(|S_tf|^2;|S_tf|>ce^{\delta t }\right) \\
  &\leq& 2t^{-1}e^{-\alpha t}\mathbb{P}_{\delta_x}\left(|\langle f,X_t\rangle-\langle T^\alpha_{\epsilon_t}f,X_{[t]}\rangle|^2 \right) +2\frac{[t]}{t}e^{-\alpha\epsilon_t}\mathbb{P}_{\delta_x}\left(|S_{[t]}(T^\alpha_{\epsilon_t}f)|^2;|S_tf|>ce^{\delta t }\right)\\
  &\leq& 2t^{-1}e^{-\alpha t}\mathbb{P}_{\delta_x}\left(|\langle f,X_t\rangle-\langle T^\alpha_{\epsilon_t}f,X_{[t]}\rangle|^2 \right)\\ &&+2\frac{[t]}{t}e^{-\alpha\epsilon_t}\mathbb{P}_{\delta_x}\left(|S_{[t]}(T^\alpha_{\epsilon_t}f)|^2;
  |S_{[t]}(T^\alpha_{\epsilon_t}f)|>ce^{\alpha\epsilon_t/2}e^{\delta [t] }\right)\\
  &&+2\frac{[t]}{t}e^{-\alpha\epsilon_t}\mathbb{P}_{\delta_x}\left(|S_{[t]}(T^\alpha_{\epsilon_t}f)|^2;
  |S_{[t]}(T^\alpha_{\epsilon_t}f)|\leq ce^{\alpha\epsilon_t/2}e^{\delta [t]},|S_tf|>ce^{\delta t}\right)\\
  &=:& A_1(t)+A_2(t)+A_3(t).
\end{eqnarray*}
 To prove \eqref{4.1} we only need to prove that  $\lim_{t\to\infty}A_j(t)=0$ for $j=1,2,3$.
In the following we give the detailed proof of $\lim_{t\to\infty}A_2(t)=0$.
The arguments for $A_1(t)$ and $A_3(t)$ are similar and are omitted. To prove $\lim_{t\to\infty}A_2(t)=0$
we only need to prove, for $m\in\mathbb{N}$,
\begin{equation}\label{4.38}
\lim_{m\to\infty}\sup_{0\leq s<1}\mathbb{P}_{\delta_x}\left(|S_m(T^\alpha_sf)|^2;|S_m(T^\alpha_sf)|>ce^{\delta m}\right)=0.
\end{equation}
Let
 $$
 F(t,f,c, \delta):=\mathbb{P}_{\delta_x}\left(|S_tf|^2;|S_tf|>ce^{\delta t }\right).
 $$
   Then \eqref{4.38} is equivalent to
 \begin{equation}\label{F(m)to0}
\lim_{m\to\infty}\sup_{0\leq s<1}F(m,T^\alpha_sf,c, \delta)=0.
 \end{equation}
Note that
\begin{eqnarray}\label{4.2}
  S_{m+1}(T^\alpha_sf)&=&\left(\frac{1}{m+1}\right)^{1/2}e^{-(\alpha/2) (m+1)}\left(\langle T^\alpha_sf, X_{m+1}\rangle-\langle T^\alpha_{s+1}f, X_m\rangle\right)\nonumber\\
  &&+\left(\frac{1}{m+1}\right)^{1/2}e^{-(\alpha/2) (m+1)}\left(\langle T^\alpha_{s+1}f, X_m\rangle-T^\alpha_{m+s+1}f(x)\right)\nonumber\\
  &=&\left(\frac{1}{m+1}\right)^{1/2}R(m,T^\alpha_sf)+\left(\frac{m}{m+1}\right)^{1/2}e^{-\alpha/2}S_m(T^\alpha_{s+1}f),
\end{eqnarray}
where $R(t,f):=e^{-(\alpha/2) (t+1)}\left(\langle f,X_{t+1}\rangle-\langle T^\alpha_1f, X_t\rangle\right)$.
Thus we have
 \begin{eqnarray*}
 &&F(m+1,T^\alpha_sf,c, \delta)\\
 &\leq& \mathbb{P}_{\delta_x}\left(|S_{m+1}(T^\alpha_sf)|^2;|S_m(T^\alpha_{s+1}f)|>ce^{\alpha/2}e^{\delta m }\right)\\
 &&+\mathbb{P}_{\delta_x}\left(|S_{m+1}(T^\alpha_sf)|^2;|S_m(T^\alpha_{s+1}f)|\leq ce^{\alpha/2}e^{\delta m},|S_{m+1}(T^\alpha_sf)|>ce^{\delta (m+1)}\right)\\
 &=:&M_1(m,T^\alpha_sf,c,\delta)+M_2(m,T^\alpha_sf,c,\delta).
 \end{eqnarray*}
 Put
 \begin{eqnarray*}
    A_1(m,f,c, \delta) &=& \{|S_m(T^\alpha_1f)|>ce^{\alpha/2}e^{\delta m }\}, \\
   A_2(m,f,c, \delta) &=& \{|S_m(T^\alpha_1f)|\leq ce^{\alpha/2}e^{\delta m }, |S_{m+1}f|> ce^{\delta(m+1)}\}.
 \end{eqnarray*}
 Since
$A_1(m,T^\alpha_sf,c, \delta)\in\mathcal{F}_m$ and $\mathbb{P}_{\delta_x}(R(m,T^\alpha_sf) |\mathcal{F}_m)$=0,
we have by \eqref{4.2} that
 \begin{eqnarray*}
   M_1(m,T^\alpha_sf, c, \delta)&=&\frac{1}{m+1}\mathbb{P}_{\delta_x} \left(|R(m,T^\alpha_sf)|^2;A_1(m,T^\alpha_sf,c,\delta)\right)\nonumber\\
&& +\frac{m}{m+1}e^{-\alpha}\mathbb{P}_{\delta_x}\left(|S_m(T^\alpha_{s+1}f)|^2;A_1(m,T^\alpha_sf,c, \delta)\right)\nonumber\\
   &=&
  \frac{1}{m+1}\mathbb{P}_{\delta_x}\left(|R(m,T^\alpha_sf)|^2;A_1(m,T^\alpha_sf,c, \delta)\right)\nonumber\\
&& +\frac{m}{m+1}e^{-\alpha}F\left(m,T^\alpha_{s+1}f,ce^{\alpha/2}, \delta\right),
 \end{eqnarray*}
 and
 \begin{eqnarray*}
   M_2(m,T^\alpha_sf, c, \delta)&\leq&\frac{2}{m+1}\mathbb{P}_{\delta_x}\left(|R(m,T^\alpha_sf)|^2;A_2(m,T^\alpha_sf,c, \delta)\right)\\
   &&+\frac{2m}{m+1}e^{-\alpha}\mathbb{P}_{\delta_x}\left(|S_m(T^\alpha_{s+1}f)|^2;A_2(m,T^\alpha_sf,c, \delta)\right).
 \end{eqnarray*}
 Thus we have
 \begin{eqnarray}\label{4.4}
   F(m+1,T^\alpha_sf,c, \delta)&\leq& \frac{m}{m+1}e^{-\alpha}F\left(m,T^\alpha_{s+1}f,ce^{\alpha/2}, \delta\right)
   \nonumber\\
   &&+\frac{1}{m+1}(G_1(m,T^\alpha_sf,c, \delta)+G_2(m,T^\alpha_sf,c, \delta)),
 \end{eqnarray}
 where
 \begin{eqnarray*}
   G_1(m,T^\alpha_sf,c, \delta) &=& 2\mathbb{P}_{\delta_x}\left(|R(m,T^\alpha_sf)|^2;A_1(m,T^\alpha_sf,c, \delta)\cup A_2(m,T^\alpha_sf,c, \delta)\right),\\
   G_2(m,T^\alpha_sf,c, \delta) &=& 2me^{-\alpha}\mathbb{P}_{\delta_x}\left(|S_m(T^\alpha_{s+1}f)|^2;A_2(m,T^\alpha_sf,c, \delta)\right).
 \end{eqnarray*}
 Iterating \eqref{4.4}, we get
 \begin{eqnarray}\label{4.5}
    F(m+1,T^\alpha_sf,c, \delta)&\leq&\frac{1}{m+1}\sum_{k=0}^me^{-k\alpha}G_1(m-k,T^\alpha_{k+s}f,ce^{\alpha k/2}, \delta)\nonumber\\
&&    +\frac{1}{m+1}\sum_{k=0}^me^{-k\alpha}G_2(m-k,T^\alpha_{k+s}f,ce^{\alpha k/2}, \delta)\nonumber\\
    &=:&L_1(s,f, m)+L_2(s,f, m).
 \end{eqnarray}
Therefore, to prove \eqref{F(m)to0}, we only need to prove that
\begin{equation}\label{Mto0}
\sup_{0\leq s<1}L_1(s, f, m)\to 0\quad \mbox{ and }\sup_{0\leq s<1}L_2(s, f, m)\to 0,\quad \mbox{ as } m\to\infty.
\end{equation}

  First, we consider $L_1(s, f, m)$. Let $\widetilde{f}(x)=f(x)-\sum_{|p|=\gamma(f)}a_p\phi_p(x)$. Denote
 $$
 A_{k,m,s}=A_1(m-k,T^\alpha_{k+s}f,ce^{\alpha k/2}, \delta)\cup A_2(m-k,T^\alpha_{k+s}f),ce^{\alpha k/2}, \delta).
 $$
 Then we have
 \begin{eqnarray*}
   L_1(s, f, m) &=& \frac{2}{m+1}\sum_{k=0}^me^{-k\alpha}\mathbb{P}_{\delta_x}(|R(m-k,T^\alpha_{k+s}f)|^2;A_{k,m,s}) \\
    &\leq&  \frac{4}{m+1}\sum_{k=0}^me^{-k\alpha}\mathbb{P}_{\delta_x}(|R(m-k,T^\alpha_{k+s}\widetilde{f})|^2)\\
    &&+\frac{4}{m+1}\sum_{k=0}^me^{-k\alpha}\mathbb{P}_{\delta_x}(|R(m-k,e^{\alpha (k+s)/2}\sum_{|p|=\gamma(f)}a_p\phi_p)|^2;A_{k,m,s})\\
    &\leq& \frac{4}{m+1}\sum_{k=0}^m e^{-k\alpha}\mathbb{P}_{\delta_x}(|R(m-k,T^\alpha_{k+s}\widetilde{f})|^2)\\
    &&+\frac{4e^{\alpha}}{m+1}\sum_{k=0}^m\mathbb{P}_{\delta_x}(|R(k,\sum_{|p|=\gamma(f)}a_p\phi_p)|^2;A_{m-k,m,s})\\
    &=:&L_{1,1}(s, f, m)+L_{1,2}(s,f, m).
 \end{eqnarray*}
 For $L_{1,1}(s, f, m)$, by the Markov property, we have
 \begin{eqnarray}\label{4.8}
   \mathbb{P}_{\delta_x}(|R(m-k,T^\alpha_{k+s}\widetilde{f})|^2&=&e^{-\alpha(m-k+1)}\mathbb{P}_{\delta_x}\left(\langle T^\alpha_{k+s}\widetilde{f},X_{m-k+1}\rangle-\langle T^\alpha_{k+s+1}\widetilde{f}, X_{m-k}\rangle\right)^2\nonumber\\
   &=&e^{-\alpha(m-k+1)}\mathbb{P}_{\delta_x}\langle \mathbb{V}ar_{\delta_\cdot}\langle T^\alpha_{s+k}\widetilde{f}, X_1 \rangle,~X_{m-k}\rangle\nonumber\\
   &=&e^{-\alpha}T_{m-k}(\mathbb{V}ar_{\delta_\cdot}\langle T^\alpha_{k+s}\widetilde{f},X_1\rangle)(x).
 \end{eqnarray}
By \eqref{2.40}, there exists $r\in\N$ such that $|T^\alpha_{k+s}\widetilde{f}(x)|=
e^{\alpha (k+s)}|T_{k+s}\widetilde{f}(x)|\lesssim e^{(\alpha/2)k}e^{-bk}(1+\|x\|^r)$
for $0\leq s<1$.
So by \eqref{2.24}, we obtain
 \begin{equation}\label{4.40}
   \mathbb{V}ar_{\delta_x}\langle T^\alpha_{k+s}\widetilde{f},X_1 \rangle=Ae^{\alpha}\int_0^1 e^{\alpha u}T_{1-u} [T_uT^\alpha_{k+s}\widetilde{f}]^2(x)du\lesssim e^{\alpha k}e^{-2bk}(1+\|x\|^{2r}),\quad s\in[0,1).
 \end{equation}
Thus $\mathbb{P}_{\delta_x}(|R(m-k,T^\alpha_{k+s}\widetilde{f})|^2\lesssim e^{\alpha k}e^{-2bk}(1+\|x\|^{2r})$. So,
\begin{equation}\label{4.7}
  L_{1,1}(s, f, m)\lesssim \frac{1}{m+1}\sum_{k=0}^\infty e^{-2bk}(1+\|x\|^{2r})\lesssim
  \frac{1}{m+1}(1+\|x\|^{2r})\to 0\quad m\to\infty .
\end{equation}
Now we consider $L_{1,2}(s, f, m)$.
Using \eqref{1.24} with $t=k$, $s=1$ and the function $f$ replaced by
$f_1:=\sum_{|p|=\gamma(f)}a_p\phi_p(x)$, we have
\begin{eqnarray*}
  R(k,\sum_{|p|=\gamma(f)}a_p\phi_p)&=&e^{-(\alpha/2) (k+1)}\left(\langle f_1,X_{k+1}\rangle-\langle T^\alpha_1f_1,X_k\rangle\right)\\
  &=&e^{-(\alpha/2) (k+1)}(\langle f_1,\widetilde{X}_{1}^k\rangle-{\bf P}_{\delta_x}(\langle f_1,\widetilde{X}_{1}^k\rangle|\mathcal{G}_k)+e^{-(\alpha/2) k}\sum_{u\in\cL_k}(Y_1^{u,k}-y^{u,k}_1)\\
  &=:&J_0(k)+J_1(k),
\end{eqnarray*}
where $Y_t^{u,k}, y^{u,k}_t$ are defined in \eqref{e:new}.
So for any $\epsilon>0$,
\begin{eqnarray}\label{4.9}
  L_{1,2}(s, f, m)&\leq&\frac{4e^\alpha}{m+1}\sum_{k<m\epsilon}\mathbb{P}_{\delta_x}(|R(k,\sum_{|p|=\gamma(f)}a_p\phi_p)|^2)+
  \frac{8e^\alpha}{m+1}\sum_{m\epsilon\leq k\leq m}\mathbb{P}_{\delta_x}(|J_0(k)|^2;A_{m-k,t,s})\nonumber\\
  &&+\frac{8e^\alpha}{m+1}\sum_{m\epsilon\leq k\leq m}\mathbb{P}_{\delta_x}(|J_1(k)|^2;A_{m-k,t,s})\nonumber\\
  &=:& (I)+(II)+(III).
\end{eqnarray}
Using arguments similar to those leading to \eqref{4.8} and \eqref{4.40}, we get
\begin{equation*}
  \mathbb{P}_{\delta_x}R(k,\sum_{|p|=\gamma(f)}a_p\phi_p)^2=e^{-\alpha}T_k( \mathbb{V}ar_{\delta_\cdot}\langle\sum_{|p|=\gamma(f)}a_p\phi_p, X_1\rangle)(x)\lesssim
  1+\|x\|^{2\gamma(f)}.
\end{equation*}
Thus
\begin{equation}\label{4.11}
  (I)\lesssim \epsilon(1+\|x\|^{2\gamma(f)}).
\end{equation}
For (II) and (III), we claim that
\begin{description}
  \item{(i)} $|J_0(k)|^2$ and $|J_1(k)|^2$ ,$k=1,2,...$ are both uniformly integrable with respect to $\P_{\delta_x}$;
  \item{(ii)} $\sup_{k>m\epsilon}\sup_{0\leq s<1}\mathbb{P}_{\delta_x}(A_{m-k,m,s})\to 0$ as $m\to\infty.$
\end{description}
Using the claims, we have
\begin{eqnarray*}
   &&\sup_{0\leq s<1}\frac{1}{m+1}\sum_{m\epsilon\leq k\leq m}\mathbb{P}_{\delta_x}(|J_0(k)|^2;A_{m-k,m,s})\\
   &\leq & \frac{1}{m+1}\sum_{m\epsilon\leq k\leq m}\mathbb{P}_{\delta_x}(|J_0(k)|^2;|J_0(k)|>M)\\
   &&+\frac{1}{m+1}\sum_{m\epsilon\leq k\leq m}\sup_{0\leq s<1}\mathbb{P}_{\delta_x}(|J_0(k)|^2;|J_0(k)|\leq M,A_{m-k,m,s}) \\
   &\leq &  \sup_{k\geq 1}\mathbb{P}_{\delta_x}(|J_0(k)|^2;|J_0(k)|>M)+M^2\sup_{k>m\epsilon}\sup_{0\leq s<1}\mathbb{P}_{\delta_x}(A_{m-k,m,s}).
\end{eqnarray*}
First letting $m\to\infty$ and then $M\to\infty$, we get
$$
\sup_{0\leq s<1}\frac{1}{m+1}\sum_{m\epsilon\leq k\leq m}\mathbb{P}_{\delta_x}(|J_0(k)|^2;A_{m-k,m,s})\to 0
\quad \mbox{ as } m\to\infty.
$$
Similarly, we also have
$$
\sup_{0\leq s<1}\frac{1}{m+1}\sum_{m\epsilon\leq k\leq m}\mathbb{P}_{\delta_x}(|J_1(k)|^2;A_{m-k,m,s})\to 0,
\quad \mbox{ as } m\to\infty.
$$
Thus, we have
\begin{equation}\label{4.16}
   \limsup_{m\to\infty}\sup_{0\leq s<1}L_{1,2}(s, f, m)\lesssim\epsilon(1+\|x\|^{2\gamma(f)}).
\end{equation}
Letting  $\epsilon\to0$, we get $\lim_{m\to\infty}\sup_{0\leq s<1}L_{1,2}(s,f,m)=0$.
Therefore, by \eqref{4.7}, we get
$$
\lim_{m\to\infty}\sup_{0\leq s<1}L_1(s,f,m)=0.
$$

Now we prove the claims (i) and (ii).

(i) For $J_0(k)$, by \eqref{4.15}, there exists $r\in \N$ such that
\begin{eqnarray*}
  \mathbb{P}_{\delta_x}|J_0(k)|^4&\leq& e^{-2\alpha(k+1)}\mathbb{P}_{\delta_x}\left(\left\langle 1+\|\cdot\|^{4r},X_k \right\rangle+\left\langle 1+\|\cdot\|^{2r},X_k \right\rangle^2\right)\\
  &\lesssim&e^{-(k+2)\alpha}(1+\|x\|^{4r})+e^{-2\alpha(k+1)}\mathbb{P}_{\delta_x}\left\langle 1+\|\cdot\|^{2r},X_k \right\rangle^2.
\end{eqnarray*}
By\eqref{e:ousmgp} and \eqref{2.24}, we get
$$
\mathbb{P}_{\delta_x}\langle(1+\|\cdot\|^{2r}),X_k\rangle=e^{\alpha k}T_k(1+\|\cdot\|^{2r})(x)\lesssim e^{\alpha k}( 1+\|x\|^{2r}),
$$
$${\mathbb V}ar_{\delta_x}\langle(1+\|\cdot\|^{2r}),X_k\rangle\lesssim e^{2\alpha k}(1+\|x\|^{4r}).$$
So we have
\begin{equation}\label{9.6}
  \mathbb{P}_{\delta_x}\langle(1+\|\cdot\|^{2r}),X_k\rangle^2\lesssim  e^{2\alpha k}(1+\|x\|^{4r}).
\end{equation}
Thus $\sup_{k>0}\mathbb{P}_{\delta_x}|J_0(k)|^4<\infty$ which implies $|J_0(k)|^2$ is uniformly integrable.

For $J_1(k)$, from the proof of \eqref{10.5}, we see that
 \eqref{10.5} is also true when $\alpha=2\gamma(f)b$.
So we have $J_1(k)\stackrel{d}{\to}\sqrt{W_\infty}G$
where $G$ is a Gaussian random variable.
We  also have $\mathbb{P}_{\delta_x}|J_1(k)|^2\to \mathbb{P}_{\delta_x}W_\infty G^2$.
Thus, $J_1(k)$ is uniformly integrable
by \cite[Theorem 5.5.2]{Durrett} and Skorokhod's representation theorem.

(ii) Recall that
\begin{eqnarray*}
  A(m-k,m,s)&=& A_1(k,T^\alpha_{m+s-k}f,ce^{\alpha (m-k)/2}, \delta)\cup A_2(k,T^\alpha_{m+s-k}f,ce^{\alpha(m-k)/2}, \delta).
\end{eqnarray*}
By Chebyshev's inequality
\begin{eqnarray*}
 \mathbb{P}_{\delta_x}( A_1(k,T^\alpha_{m+s-k}f,ce^{\alpha (m-k)/2}), \delta) &\leq& c^{-2}
 e^{-\alpha(m-k+1)}e^{-2\delta k }\mathbb{P}_{\delta_x}|S_k(T^\alpha_{m+s-k+1}f)|^2.
\end{eqnarray*}
By \eqref{2.24} and \eqref{2.40}, we have
\begin{eqnarray}\label{8.51}
  &&\mathbb{P}_{\delta_x}|S_k(T^\alpha_{m+s-k+1}f)|^2=k^{-1}e^{-\alpha k}\mathbb{V}ar_{\delta_x}\langle T^\alpha_{m+s-k+1}f,X_k \rangle\nonumber\\
  &=&Ak^{-1}e^{2\alpha(m+s-k+1)}\int_0^k e^{\alpha u}T_{k-u} [T_{u+m+s-k+1}f]^2(x)ds\lesssim e^{\alpha(m+s-k+1)}(1+\|x\|^{2r}),
\end{eqnarray}
which implies
\begin{equation}\label{4.22}
   \sup_{k>m\epsilon}\sup_{0\leq s<1}\mathbb{P}_{\delta_x}( A_1(k,T^\alpha_{m+s-k}f,ce^{\alpha (m-k)/2}), \delta)
  \lesssim \sup_{k>m\epsilon}e^{-2\delta k}(1+\|x\|^{2r})\to 0,\quad m\to\infty.
\end{equation}
It is easy to see that
\begin{equation}\label{8.50}
 A_2(k,T^\alpha_{m+s-k}f,ce^{\alpha(m-k)/2}, \delta)\subset\left\{|R(k,T^\alpha_{m+s-k}{f})|>ce^{\alpha(m-k)/2}
e^{\delta k}(e^{\delta}\sqrt{k+1}-\sqrt{k})\right\}.
\end{equation}
Similarly, by Chebyshev's inequality, we have
\begin{equation*}
   \mathbb{P}_{\delta_x}(A_2(k,T^\alpha_{m+s-k}f,ce^{\alpha (m-k)/2}), \delta)\leq c^{-2}e^{-\alpha(m-k)}
  e^{-2\delta k}(e^{\delta}\sqrt{k+1}-\sqrt{k})^{-2}\mathbb{P}_{\delta_x}|R(k,T^\alpha_{m+s-k}{f})|^2.
\end{equation*}
Using an argument similar to that leading to \eqref{4.8}, we get
$$
\mathbb{P}_{\delta_x}|R(k,T^\alpha_{m+s-k}{f})|^2=e^{-\alpha}T_{k}(\mathbb{V}ar_{\delta_\cdot}\langle T^\alpha_{m+s-k}f,X_1\rangle)(x),
$$
so using an argument similar to that leading to \eqref{4.40}, we obtain
\begin{equation}\label{4.26}
  \mathbb{V}ar_{\delta_x}\langle T^\alpha_{m+s-k}f,X_1\rangle\lesssim e^{\alpha(m-k)}(1+\|x\|^{2r}),
\end{equation}
which implies $\mathbb{P}_{\delta_x}|R(k,T^\alpha_{m+s-k}{f})|^2\lesssim e^{\alpha(m-k)}(1+\|x\|^{2r}).$
Thus
\begin{equation}\label{4.27}
  \sup_{k>t\epsilon}\sup_{0\leq s<1}\mathbb{P}_{\delta_x}(A_2(k,T^\alpha_{m+s-k}f,ce^{\alpha(m-k)/2}))\lesssim \sup_{k>m\epsilon}e^{-2\delta k}(e^{\delta}\sqrt{k+1}-\sqrt{k})^{-2}(1+\|x\|^{2r})\to0,
  \end{equation}
as $m\to\infty.$
Claim (ii) now follows easily from \eqref{4.22} and \eqref{4.27}.

To finish the proof, we need to show that
\begin{equation}\label{4.25}
  \sup_{0\leq s<1}L_2(s,f, m)=\sup_{0\leq s<1}\frac{1}{m+1}\sum_{k=0}^me^{-\alpha(m-k)}
   G_2(k,T^\alpha_{m+s-k}f,ce^{\alpha (m-k)/2}, \delta)\to0,\quad m\to\infty.
\end{equation}
By \eqref{8.50} and Chebyshev's inequality, we have
\begin{eqnarray*}
  &&G_2(k,T^\alpha_{m+s-k}f,ce^{\alpha (m-k)/2}, \delta)\\
  &=&2e^{-\alpha}k\mathbb{P}_{\delta_x}\left(|S_k(T^\alpha_{m+s-k+1}f)|^2;
     A_2(k,T^\alpha_{m+s-k}f,ce^{\alpha(m-k)/2}, \delta)\right)\\
  &\leq& 2e^{-\alpha}kce^{\alpha(m-k+1)/2}e^{\delta k }\mathbb{P}_{\delta_x}\left(|S_k(T^\alpha_{m+s-k+1}f)|;
    A_2(k,T^\alpha_{m+s-k}f,ce^{\alpha(m-k)/2}, \delta)\right)\\
%  &\leq &2e^{-\alpha}kce^{\alpha/2(m-k+1)}e^{\delta k }\mathbb{P}_{\delta_x}\left(|S_k(T^\alpha_{m+s-k+1}f)|;
%  |R(k,T^\alpha_{m+s-k}{f})|>ce^{\alpha/2(m-k)}e^{\delta k}
%  (e^{\delta}\sqrt{k+1}-\sqrt{k})\right)\\
  &\leq &2kce^{\delta k+\alpha(m-k-1)/2}\mathbb{P}_{\delta_x}\left(|S_k(T^\alpha_{m+s-k+1}f)|;
  |R(k,T^\alpha_{m+s-k}{f})|>ce^{\delta k+\alpha/2(m-k)}
  (e^{\delta}\sqrt{k+1}-\sqrt{k})\right)\\
  &\leq& 2c^{-1}e^{-\alpha(m-k+1)/2-\delta k} (e^{\delta}\sqrt{k+1}-\sqrt{k})^{-2}k\mathbb{P}_{\delta_x}\left(|S_k(T^\alpha_{m+s-k+1}f)|
  |R(k,T^\alpha_{m+s-k}f)|^2\right)\\
  &\lesssim &e^{-\alpha(m-k)/2-\delta k}\mathbb{P}_{\delta_x}\left(|S_k(T^\alpha_{m+s-k+1}f)||R(k,T^\alpha_{m+s-k}f)|^2\right)\\
  &=& e^{-\alpha(m-k)/2-\delta k}e^{-\alpha(k+1)}\mathbb{P}_{\delta_x}\left(|S_k(T^\alpha_{m+s-k+1}f)|\langle\mathbb{V}ar_{\delta_\cdot}\langle T^\alpha_{m+s-k}f,X_1\rangle,X_k\rangle\right).
\end{eqnarray*}
By \eqref{4.26}, we get
\begin{eqnarray*}
   &&\mathbb{P}_{\delta_x}\left(|S_k(T^\alpha_{m+s-k+1}f)|\langle\mathbb{V}ar_{\delta_\cdot}\langle T^\alpha_{m+s-k}f,X_1\rangle,X_k\rangle\right)\\
   &\lesssim& e^{\alpha(m-k)}\mathbb{P}_{\delta_x}\left(|S_k(T^\alpha_{m+s-k+1}f)|\langle(1+\|\cdot\|^{2r}),X_k\rangle\right) \\
   &\leq &  e^{\alpha(m-k)}\sqrt{\mathbb{P}_{\delta_x}\left(|S_k(T^\alpha_{m+s-k+1}f)|^2\right)
   \mathbb{P}_{\delta_x}\langle(1+\|\cdot\|^{2r}),X_k\rangle^2}.
\end{eqnarray*}
Thus by \eqref{8.51} and  \eqref{9.6}, we get
\begin{equation*}
  \mathbb{P}_{\delta_x}\left(|S_k(T^\alpha_{m+s-k+1}f)|\langle\mathbb{V}ar_{\delta_\cdot}\langle T^\alpha_{m+s-k}f,X_1\rangle,X_k\rangle\right)\lesssim e^{\alpha(m-k)}e^{\alpha(m+k)/2}(1+\|x\|^{3r}),
\end{equation*}
which implies
$$
G_2(k,T^\alpha_{m+s-k}f,ce^{\alpha (m-k)/2})\lesssim e^{\alpha(m-k)}e^{-\delta k}(1+\|x\|^{3r}).
$$
Therefore, we have
\begin{equation}\label{4.34}
   \sup_{0\leq s<1}L_2(s,f, m)\lesssim\frac{1}{m+1}\sum_{k=0}^me^{-\delta k}(1+\|x\|^{3r})\to 0,\quad m\to\infty.
\end{equation}
Hence, $\lim_{t\to\infty}A_2(t)=0$. \hfill$\Box$

\smallskip

In the following lemma we give a result similar to Lemma \ref{lem:5.5} for the process $I$.

\begin{lemma}\label{lem:5.6}
Assume $f\in \mathcal{P}$ satisfies $\alpha = 2\gamma(f)b$. Define
$$
Y^*_t(f):=t^{-1/2}e^{-(\alpha/2) t}\left(\langle f,I_t\rangle-\mathbb{Q}_{\delta_x}\langle f,I_t\rangle\right).
$$
Then for any $c>0$ and $\delta>0$, we have
\begin{equation}\label{4.41}
\lim_{t\to\infty}\mathbb{Q}_{\delta_x}\left(|Y^*_t(f)|^2;|Y^*_t(f)|>ce^{\delta t }\right)=0.
\end{equation}
\end{lemma}

\textbf{Proof:}\quad
Recall the decomposition in \eqref{12}.
Define
$$
S^*_t=t^{-1/2}e^{-(\alpha/2) t}(\langle f,\widetilde{X}_t\rangle-{\bf P}_{\delta_x}\langle f,\widetilde{X}_t\rangle),
$$
$$
S_t=t^{-1/2}e^{-(\alpha/2) t}(\langle f,\Lambda_t\rangle-{\bf P}_{\delta_x}\langle f,\Lambda_t\rangle),
$$
and
$$
\widetilde{Y}_t=t^{-1/2}e^{-(\alpha/2) t}(\langle f,I_t\rangle-{\bf P}_{\delta_x}\langle f,I_t\rangle).
$$
Then we have $\widetilde{Y}_t=S_t-S_t^*$. Thus
\begin{eqnarray*}\label{4.42}
  {\bf P}_{\delta_x}(|\widetilde{Y}_t|^2;|\widetilde{Y}_t|>ce^{\delta t})&\leq& {\bf P}_{\delta_x}(|\widetilde{Y}_t|^2;|S_t|>(c/2)e^{\delta t})+{\bf P}_{\delta_x}(|\widetilde{Y}_t|^2;|S^*_t|>(c/2)e^{\delta t})\\
  &\leq& 2{\bf P}_{\delta_x}(|S_t|^2;|S_t|>(c/2)e^{\delta t})+2{\bf P}_{\delta_x}(|S^*_t|^2)+{\bf P}_{\delta_x}(|\widetilde{Y}_t|^2;|S^*_t|>(c/2)e^{\delta t})\\
  &=& I_1(t)+I_2(t)+I_3(t).
\end{eqnarray*}
By Lemma \ref{lem:5.5},
we have $\lim_{t\to\infty}I_1(t)=0$.
By \eqref{1.8}, we have
\begin{equation*}
I_2(t)=2t^{-1}e^{-\alpha t}{\bf V}ar_{\delta_x}\langle f, \widetilde{X}_t\rangle\to0,\quad t\to\infty.
\end{equation*}
Since
$I_t$ and $\widetilde{X}$ are independent, we have
\begin{equation*}
  I_3(t)={\bf P}_{\delta_x}(|\widetilde{Y}_t|^2)~{\bf P}_{\delta_x}(|S^*_t|>(c/2)e^{\delta t}).
\end{equation*}
Since $S_t=S^*_t+\widetilde{Y}_t$, and $S^*_t$ and $\widetilde{Y}_t$ are independent, by \eqref{3.11}, we get
\begin{equation*}
  {\bf P}_{\delta_x}(|\widetilde{Y}_t|^2)={\bf P}_{\delta_x}(|S_t|^2)-{\bf P}_{\delta_x}(S^*_t|^2)\to \rho_f^2,\quad t\to\infty.
\end{equation*}
By Chebyshev's inequality, we have
$$
{\bf P}_{\delta_x}(|S^*_t|>(c/2)e^{\delta t})\leq (c/2)^{-2}e^{-2\delta t}){\bf P}_{\delta_x}(|S^*_t|^2)\to0,\quad t\to\infty.
$$
Hence $\lim_{t\to\infty}I_3(t)=0$. Thus
\begin{equation}\label{4.45}
  {\bf P}_{\delta_x}(|\widetilde{Y}_t|^2;|\widetilde{Y}_t|>ce^{\delta t})\to 0.
\end{equation}

Recall that under ${\bf P}_{\delta_x}$, $I_t=\sum_{j=1}^{N}I^{j}_t,$  where $I^{j},j=1,...$
are independent copies of $I$ under $\mathbb{Q}_{\delta_x}$, and are independent of $N$. Thus,
\begin{equation*}
  {\bf P}_{\delta_x}(|\widetilde{Y}_t|^2;|\widetilde{Y}_t|>ce^{\delta t})\geq {\bf P}_{\delta_x}(|\widetilde{Y}_t|^2;|\widetilde{Y}_t|>ce^{\delta t},N=1)={\bf P}_{\delta_x}(N=1)\mathbb{Q}_{\delta_x}\left(|Y^*_t(f)|^2;|Y^*_t(f)|>ce^{\delta t }\right).
\end{equation*}
Since ${\bf P}_{\delta_x}(N=1)>0$,
\eqref{4.41} follows easily from \eqref{4.45}. \hfill$\Box$
\smallskip

Now, we are ready to prove Theorem \ref{The:1.5}.

\textbf{Proof of Theorem \ref{The:1.5}:}\quad
The proof is similar to that of Theorem \ref{The:1.4}.
 We define an ${\mathbb R}^2$-valued random variable by
\begin{equation*}
  U_1(t):=(e^{-\alpha (t)}\left\|\Lambda_{t}\right\|, t^{-1/2}e^{-(\alpha/2) (t)}\langle f,\Lambda_{t}\rangle).
\end{equation*}
We need to show that as $t\to\infty$,
\begin{equation}\label{7.5}
  U_1(t)\stackrel{d}{\to}(W_\infty, \sqrt{W_\infty}G_2(f)),
  \end{equation}
 where $G_2(f)\sim\mathcal{N}(0,\rho_f^2)$.
Let $n>0$ and write
\begin{equation*}
  U_1(nt)=(e^{-\alpha (nt)}\left\|\Lambda_{nt}\right\|, (nt)^{-1/2}e^{-(\alpha/2) (nt)}\langle f,\Lambda_{nt}\rangle).
\end{equation*}
Recall the representation \eqref{3.6}.
Define
$$
Y_t^{u, n}:=((n-1)t)^{-1/2}e^{-\alpha(n-1)t/2}\langle f,I_{(n-1)t}^{u,t}\rangle\quad \mbox{and}
\quad y_t^{u, n}:={\bf P}_{\mu}(Y^{u, n}_t|\mathcal{G}_t).
$$
$Y_t^{u, n}$ has the same distribution as $Y^n_t:=((n-1)t)^{-1/2}e^{-\alpha(n-1)t/2}\langle f,I_{(n-1)t}\rangle$
under $\mathbb{Q}_{\delta_{Z_u(t)}}$. Thus
\begin{eqnarray}\label{9.1}
     &&(nt)^{-1/2}e^{-(\alpha/2)nt}\langle f,\Lambda_{nt}\rangle\nonumber\\
  &=& (nt)^{-1/2}e^{-(\alpha/2)nt} \langle f,\widetilde{X}_{(n-1)t}^t\rangle+\sqrt{\frac{n-1}{n}}
    e^{-(\alpha/2) t}\sum_{u\in\cL_t}Y_t^{u, n}\nonumber\\
  &=& (nt)^{-1/2}e^{-(\alpha/2)nt}(\langle f,\widetilde{X}_{(n-1)t}^t\rangle-{\bf P}_{\mu}(\langle f,\widetilde{X}_{(n-1)t}^t\rangle|\mathcal{G}_t)+\sqrt{\frac{n-1}{n}}
   e^{-(\alpha/2) t}\sum_{u\in\cL_t}(Y_t^{u, n}-y^{u, n}_t)\nonumber\\
  &&+ (nt)^{-1/2}e^{-(\alpha/2)nt}{\bf P}_{\mu}(\langle f,\Lambda_{nt}\rangle|\mathcal{G}_t)\nonumber\\
  &=:& J_0^n(t)+J_1^n(t)+J_2^n(t).
\end{eqnarray}
Put
$\widetilde{V}_s(x):={\bf V}ar_{\delta_x}\langle f,\widetilde{X}_{s}\rangle$.
Then by \eqref{1.8},
there exists $r\in\N$ such that $\widetilde{V}_s(x)\lesssim e^{-\alpha^*s}(1+\|x\|^{2r})$.
From the definition of $\widetilde{X}^t_{s}$, we have
\begin{eqnarray}
  {\bf P}_\mu J_0^n(t)^2&=&(nt)^{-1}e^{-\alpha(nt)}{\bf P}_\mu (\langle \widetilde{V}_{(n-1)t},\Lambda_t\rangle)=(nt)^{-1}e^{-\alpha (n-1)t}\langle T_t(\widetilde{V}_{(n-1)t}, \mu\rangle\nonumber\\
  &\lesssim& (nt)^{-1}e^{-\alpha (n-1)t}e^{-\alpha^*(n-1)t}\to 0,\quad \mbox{as } t\to\infty.\label{9.3}
\end{eqnarray}
Using an argument similar to that in the proof of  Theorem \ref{The:1.4}, we can get
\begin{eqnarray}\label{9.2}
  {\bf P}_{\mu}J_{2}^n(t)^2 &=& A (nt)^{-1}e^{\alpha (n-1)t}\int_{\mathbb{R}^d}\int_0^t e^{\alpha u}T_{t-u}[T_{(n-1)t+u}f]^2(x)\,ds\,\mu(dx)\nonumber\\
  &\lesssim& n^{-1}.
\end{eqnarray}
Combining \eqref{9.3} and \eqref{9.2}, there exists $c>0$ such that
\begin{equation}\label{9.4}
 \limsup_{t\to\infty}{\bf P}_{\mu}(J_0^n(t)+J_2^n(t))^2\leq c/n.
\end{equation}

Now we consider $J_1^n(t)$.
We define an ${\mathbb R}^2$-valued random variable $U_2(n,t)$ by
\begin{eqnarray*}
U_2(n,t):=\left(e^{-\alpha t}\|\Lambda_t\|,e^{-(\alpha/2) t}\sum_{u\in\cL_t}(Y_t^{u, n}-y^{u, n}_t) \right).
\end{eqnarray*}
We claim that
\begin{equation}\label{9.5}
  U_2(n,t)\stackrel{d}{\to}(W_\infty, \sqrt{W_\infty}G_2(f)), \quad \mbox{ as } t\to\infty.
\end{equation}
Denote the characteristic function of $U_2(n,t)$ under ${\bf P}_\mu$ by
$\kappa_2(\theta_1,\theta_2,n,t)$.
Using an argument similar to that leading to \eqref{10.9}, we get
\begin{eqnarray*}
  \kappa_2(\theta_1,\theta_2,n,t)
  &=& {\bf P}_{\mu}\left(\exp\{i\theta_1e^{-\alpha t}\|\Lambda_t\|\}\exp\left\{\lambda^*\langle
      h^n_t(\cdot,~e^{-(\alpha/2)t}\theta_2)-1,\Lambda_t\rangle\right\}\right),
\end{eqnarray*}
where $h^n_t(x,\theta)=\mathbb{Q}_{\delta_x}e^{i\theta(Y^n_t-\mathbb{Q}_{\delta_x}Y^n_t)}$.
Define
\begin{equation*}
   e^n_t(x,\theta):=h^n_t(x,\theta)-1+\frac{1}{2}\theta^2\mathbb{V}_{\delta_x}Y^n_t
\end{equation*}
and $V^n_t(x):=\mathbb{V}_{\delta_x}Y^n_t$. Then
\begin{eqnarray*}
   &&\exp\left\{\lambda^*\langle h^n_t(\cdot,~e^{-(\alpha/2)t}\theta_2)-1,\Lambda_t\rangle\right\}\\
   &=& \exp\left\{-\frac{1}{2}\lambda^*\theta_2^2e^{-\alpha t}\langle V^n_t, \Lambda_t\rangle\right\}
    \exp\left\{\lambda^*\langle e^n_t(\cdot,e^{-(\alpha/2)t}\theta_2), \Lambda_t\rangle\right\}\\
   &=:&J_{1,1}(n,t)J_{1,2}(n,t).
\end{eqnarray*}

We first consider $J_{1,1}(n,t)$.
By \eqref{1.17}, we have that as $t\to\infty$,
\begin{equation*}
   e^{-\alpha t}\langle |\lambda^*V^n_t-\rho_f^2|, \Lambda_t\rangle\lesssim  t^{-1}
  e^{-\alpha t}\langle (1+\|x\|^{r}), \Lambda_t\rangle\to 0 \quad \mbox{in probability}.
\end{equation*}
It follows that
\begin{equation}\label{8.15}
  \lim_{t\to\infty}e^{-\alpha t}\langle \lambda^*V^n_t, \Lambda_t\rangle=\lim_{t\to\infty}
  e^{-\alpha t}\langle \rho_f^2, \Lambda_t\rangle=\rho_f^2W_\infty \quad \mbox{in probability},
\end{equation}
which implies that $\lim_{t\to\infty}J_{1,1}(n,t)=\exp\left\{-\frac{1}{2}\theta_2^2\rho_f^2W_\infty\right\}$.

For $J_{1,2}(n,t)$, by \eqref{3.20}, we have, for any $\epsilon>0$,
\begin{eqnarray*}
|e^n_t(x,e^{-(\alpha/2)t}\theta_2)|&\leq& \frac{1}{6}|\theta_2|^3e^{-\frac{3}{2}\alpha t}
 \mathbb{Q}_{\delta_x}\left(|Y^n_t-\mathbb{Q}_{\delta_x}Y^n_t|^3;
|Y^n_t-\mathbb{Q}_{\delta_x}Y^n_t|<\epsilon e^{\alpha t/2}\right)\\
&&+\theta_2^2e^{-\alpha t}\mathbb{Q}_{\delta_x}
\left(|Y^n_t-\mathbb{Q}_{\delta_x}Y^n_t|^2; |Y^n_t-\mathbb{Q}_{\delta_x}Y^n_t|\geq\epsilon e^{\alpha t/2}\right)\\
&\leq&  \frac{\epsilon}{6}|\theta|_2^3e^{-\alpha t}
\mathbb{Q}_{\delta_x}\left(|Y^n_t-\mathbb{Q}_{\delta_x}Y^n_t|^2\right)\\
&&+\theta_2^2e^{-\alpha t}\mathbb{Q}_{\delta_x}
\left(|Y^n_t-\mathbb{Q}_{\delta_x}Y^n_t|^2; |Y^n_t-\mathbb{Q}_{\delta_x}Y^n_t|\geq\epsilon e^{\alpha t/2}\right)\\
&=&\frac{\epsilon}{6}|\theta|_2^3e^{-\alpha t}V^n_t(x)+\theta_2^2e^{-\alpha t}F_t^n(x),
\end{eqnarray*}
where $F^n_t(x)=\mathbb{Q}_{\delta_x}\left(|Y^n_t-\mathbb{Q}_{\delta_x}Y^n_t|^2; |Y^n_t-\mathbb{Q}_{\delta_x}Y^n_t|
\geq\epsilon e^{\alpha t/2}\right)$.
It follows from Lemma \ref{lem:5.6} that $\lim_{t\to\infty}F^n_t(x)= 0$.
By \eqref{1.19}, we also have
$$
F^n_t(x)\leq\mathbb{Q}_{\delta_x}\left(|Y^n_t-\mathbb{Q}_{\delta_x}Y^n_t|^2\right) \lesssim 1+\|x\|^{2r}.
$$
Note that
$$
e^{-\alpha t} {\bf P}_{\mu}\langle F_t^n(x), \Lambda_t\rangle=\langle T_t(F^n_t),\mu\rangle.
$$
Thus by the dominated convergence theorem,
we get $\lim_{t\to\infty}e^{-\alpha t} {\bf P}_{\mu}\langle F_t^n(x), \Lambda_t\rangle= 0$.
It follows that $e^{-\alpha t}\langle F_t^n(x), \Lambda_t\rangle\to 0$ in probability.
Furthermore from \eqref{8.15}, we obtain that as $t\to\infty$,
$$
\frac{\epsilon}{6}\theta_2^3e^{-\alpha t}\langle V^n_t,
\Lambda_t\rangle\to\frac{\epsilon}{6\lambda^*}\theta_2^3 \rho_f^2W_\infty\quad \mbox{in probability}.
$$
Thus, letting $\epsilon\to 0$, we get that as $t\to\infty$,
\begin{equation}\label{4.48}
 \langle |e^n_t(x,e^{-(\alpha/2)t}\theta_2)|,\Lambda_t\rangle\to0 \quad \mbox{in probability},
\end{equation}
which implies $J_{1,2}(n,t)\to 1$ in probability, as $t\to\infty$.

Thus, when $t\to\infty$,
\begin{equation}\label{9.8}
  \exp\left\{\lambda^*\langle h^n_t(\cdot,~e^{-(\alpha/2)t}\theta_2)-1,\Lambda_t\rangle\right\}
  \to\exp\left\{-\frac{1}{2}\theta_2^2\rho_f^2 W_\infty\right\}
\end{equation}
 in probability. Since $h^n_t(x,\theta)$ is a characteristic function,
its real part is less than 1, which implies
$$
|\exp\left\{\lambda^*\langle h^n_t(\cdot,~e^{-(\alpha/2)t}\theta_2)-1,\Lambda_t\rangle\right\}|\leq 1.
$$
So by the dominated convergence theorem, we get that
\begin{equation}\label{9.10}
\lim_{t\to\infty}\kappa_2(\theta_1,\theta_2,n,t)={\bf P}_{\mu}\exp\left\{i\theta_1W_\infty\right\}
  \exp\left\{-\frac{1}{2}\theta_2^2\rho_f^2  W_\infty\right\},
 \end{equation}
which implies our claim \eqref{9.5}.
By \eqref{9.5}, we easily get that  as $t\to\infty$,
\begin{eqnarray*}
U_3(n,t)&:=&\left(e^{-\alpha (nt)}\|\Lambda_{nt}\|,J^n_1(t) \right)\stackrel{d}{\to}
(W_\infty, \sqrt{\frac{n-1}{n}}\sqrt{W_\infty}G_2(f)).
\end{eqnarray*}

Let $\mathcal{L}(nt)$ and $\widetilde{\mathcal{L}}^n(t)$ be the distributions of
$U_1(nt)$ and $U_3(n,t)$
respectively, and let $\mathcal{L}^n$ and $\mathcal{L}$ be the distributions of $(W_\infty, \sqrt{\frac{n-1}{n}}\sqrt{W_\infty}G_2(f))$ and $(W_\infty, \sqrt{W_\infty}G_2(f))$ respectively. Then, using \eqref{5.20}, we have
\begin{eqnarray}\label{4.12}
  \limsup_{t\to\infty}\beta(\mathcal{L}(nt),\mathcal{L})&\leq&
  \limsup_{t\to\infty}[\beta(\mathcal{L}(nt),\widetilde{\mathcal{L}}^n(t))+\beta(\widetilde{\mathcal{L}}^n(t),\mathcal{L}^n)+\beta(\mathcal{L}^n,\mathcal{L})]\nonumber\\
 &\leq &\limsup_{t\to\infty}({\bf P}_{\mu}(J_0^n(t)+J_2^n(t))^2)^{1/2}+0+\beta(\mathcal{L}^n,\mathcal{L}).
\end{eqnarray}
Using this and the definition of $\limsup_{t\to\infty}$, we easily get that
$$
\limsup_{t\to\infty}\beta(\mathcal{L}(t),\mathcal{L})=
\limsup_{t\to\infty}\beta(\mathcal{L}(nt),\mathcal{L})
\le \sqrt{c/n}+\beta(\mathcal{L}^n,\mathcal{L}).
$$
Letting $n\to\infty$, we get $ \limsup_{t\to\infty}\beta(\mathcal{L}(t),\mathcal{L})=0$.
The proof is now complete. \hfill$\Box$

\smallskip

\textbf{Proof of Theorem \ref{The:2.3}:}\quad First note that
\begin{eqnarray*}
  &&t^{-1/2}\|X_t\|^{-1/2}\left(\langle f,X_t\rangle-\sum_{\gamma(f)\le m<\alpha/2b}e^{(\alpha-mb)t}\sum_{|p|=m}a_pH^p_\infty\right)\\
   &=&  t^{-1/2}\|X_t\|^{-1/2}\langle f_{(cl)},X_t\rangle+ t^{-1/2}\|X_t\|^{-1/2}\left(\langle f_{(s)}, X_t\rangle-\sum_{n=1}^ke^{(\alpha-mb)t}\sum_{|p|=m}a_pH^p_\infty\right) \\
   &=:& J_1(t)+J_2(t),
\end{eqnarray*}
where $f_{(cl)}=f_{(l)}+f_{(c)}$.
By the definition of $f_{(s)}$, we have $(f_{(s)})_{(c)}=0$. Then using Theorem \ref{The:2.1} for $f_{(s)}$, we have
\begin{equation}\label{7.2}
  \|X_t\|^{-1/2}\left(\langle f_{(s)}, X_t\rangle-\sum_{n=1}^ke^{(\alpha-mb)t}\sum_{|p|=m}a_pH^p_\infty\right)\stackrel{d}{\to}G_1(f_{(s)}).
\end{equation}
 Thus
 \begin{equation}\label{7.3}
 J_2(t)\stackrel{d}{\to}0,\quad t\to\infty.
 \end{equation}
Since $\alpha=2\gamma(f_{(cl)})b$, so using Theorem \ref{The:1.5} for $f_{(cl)}$, we have
\begin{equation}\label{7.1}
(e^{-\alpha t}\|X_t\|, J_1(t)))\stackrel{d}{\to}(W^*, G_2(f_{(cl)})),
\end{equation}
where $G_2(f_{(cl)})\sim\mathcal{N}(0,\rho_{f_{(cl)}}^2)$. By \eqref{1.3}, we have
$\rho_{f_{(cl)}}^2=A\sum_{|p|=\alpha/2b}(a_p)^2$.
Combing \eqref{7.3} and \eqref{7.1}, we arrive at the conclusion of Theorem \ref{The:2.3}.\hfill$\Box$
\vspace{.1in}
\begin{singlespace}

\end{singlespace}
\end{doublespace}

\vskip 0.2truein

\noindent{\bf Yan-Xia Ren:} LMAM School of Mathematical Sciences \& Center for
Statistical Science, Peking
University,  Beijing, 100871, P.R. China. Email: {\texttt
yxren@math.pku.edu.cn}

\smallskip
\noindent {\bf Renming Song:} Department of Mathematics,
University of Illinois,
Urbana, IL 61801, U.S.A.
Email: {\texttt rsong@math.uiuc.edu}

\smallskip
 
\noindent{\bf Rui Zhang:} LMAM School of Mathematical Sciences, Peking
University,  Beijing, 100871, P.R. China. Email: {\texttt
ruizhang8197@gmail.com}
\end{document}